\RequirePackage{lineno}
\documentclass[12pt,a4paper]{article}

\usepackage[left=2cm,right=2cm, top=2cm,bottom=2cm,bindingoffset=0cm]{geometry}

\usepackage{color}

\usepackage{subfig}
\usepackage{flushend}
\usepackage{indentfirst}
\usepackage{graphics}
\usepackage{amsmath}
\usepackage{graphicx}
\usepackage{epstopdf}
\usepackage{float}
\usepackage{amssymb}
\usepackage[font=footnotesize,labelfont=bf]{caption}
\usepackage[labelsep=period]{caption}

\usepackage[font=footnotesize,labelfont=bf]{caption}
\usepackage[labelsep=period]{caption}

\graphicspath{{./figure/}}

\newtheorem{thm}{Theorem}[section]
\newtheorem{lem}{Lemma}[section]
\newtheorem{ozn}{Definition}
\newtheorem{nas}{Corollary}[section]

\newcommand{\me}{\mathbf}
\newcommand{\ld}{\left}
\newcommand{\rd}{\right}
\newcommand{\md}{\mathcal}
\newcommand{\ip}{\int_{-\pi}^{\pi}}
\begin{document}

\title{{Interpolation of functionals of stochastic sequences with stationary increments from observations with noise}}

\date{}

\maketitle

\noindent Applied Statistics. Actuarial and Financial Mathematics.  No. 2, 131--148,  2012\\
ISSN 1810-3022

\vspace{20pt}

\author{\textbf{Maksym Luz}, \textbf{Mykhailo Moklyachuk}$^{*}$    \\\\
\emph{ Department of Probability Theory, Statistics and Actuarial
Mathematics, \\
Taras Shevchenko National University of Kyiv, Kyiv 01601, Ukraine}\\
Corresponding email: moklyachuk@gmail.com}\\\\\\

\noindent \textbf{\large{Abstract}} \hspace{2pt}
The problem of optimal estimation of linear functional ${{A}_{N}}\xi =\sum\limits_{k=0}^{N}{a(k)\xi (k)}\,$ depending on the unknown values of a stochastic sequence $\xi (m)$ with stationary $n$-th increments from observations of the sequence $\xi (k)$ at points $k=-1,-2,\ldots $ and of the sequence $\xi (k)+\eta (k)$ at points of time $k=N+1,N+2,\ldots $ is considered. Formulas for calculating the mean square error and the spectral characteristic of the optimal linear estimate of the functional are proposed under condition of spectral certainty, where spectral densities of the sequences $\xi (m)$ and $\eta (m)$ are exactly known. Minimax (robust) method of estimation is applied in the case where the spectral densities are not known exactly while some sets of admissible spectral densities are given. Formulas that determine the least favorable spectral densities and the minimax spectral characteristics are proposed for some specific sets of admissible densities
\\

\noindent\textbf{{Keywords}} \hspace{2pt}
sequence with stationary increments, minimax-robust estimation, mean square error, least favorable spectral density, minimax-robust spectral characteristic.
\\

\noindent\textbf{ AMS 2010 subject classifications.} Primary: 60G10, 60G25, 60G35 Secondary: 62M20, 93E10

\noindent\hrulefill

\section{{Introduction}}

 Stochastic sequences with stationary increments are an important generalization of wide sense stationary random sequences.
 The simplest examples of sequences with stationary increments are ARIMA sequences, periodic (seasonal) time series, etc.
 Random processes with stationary increments were studied by A.M. Yaglom in the article [1].
 The article presented the spectral representation of stochastic increments, the canonical representation of the spectral density, and solution of the problem of prediction the value of the stochastic increment from known observations.
 The problem of prediction of processes with stationary increments is reduced to the problem of prediction the value of the stationary increment.
 The results of studies of random processes with stationary increments can be found in the works of M.S. Pinsker [2], A.M. Yaglom, and M.S. Pinsker [3].
 In the work [4] A.M. Yaglom used the theory of generalized random processes to study processes with stationary increments.
The classical methods of solution the problems of extrapolation, interpolation and filtering of stationary processes and sequences with known spectral densities were proposed by
A.M. Kolmogorov [5], N. Wiener [6], A.M. Yaglom [7,8]. In the article by M.M. Luz [9], the problem of extrapolation of linear functionals from unknown values of a stationary process based on known observations at discrete moments of time is solved.
\\
In the case where the exact form of the spectral density of the process is unknown, while a set of admissible spectral densities is given, the minimax method of solution the problems of optimal estimation of unknown values of a stationary process is applied.
The idea of the mentioned approach is that an estimate is determined that minimizes the value of the mean square error simultaneously for all densities from a set of admissible spectral densities.
The minimax approach to the problem of extrapolation of stationary processes was first proposed by U. Grenander [10].
In the works of Y. Franke [11], Y. Franke and H. Poor [12], S. Kassam and H. Poor [13], the problems of minimax extrapolation and interpolation of stationary sequences are solved by convex optimization methods.
In the works of M.P. Moklyachuk [14] -- [18], the problems of extrapolation, interpolation and filtration for stationary processes and sequences are investigated.
Vector stationary processes were investigated by Y.A. Rozanov [19]. In their paper O. Yu. Masyutka and M.P. Moklyachuk [20] obtained minimax estimates of linear functionals from unknown values of vector stationary processes and sequences. In the works by I.I. Dubovetska, M.P. Moklyachuk [21] and I.I. Dubovetska, O. Yu. Masyutka and M.P. Moklyachuk [22]  the problems of minimax filtering and interpolation of periodically correlated sequences are  solved. In the work of M.M. Luz and M.P. Moklyachuk [23], the problem of minimax interpolation for stochastic sequences with stationary increments for known values of this sequence is studied.

This paper presented results of investigation the problem of optimal linear interpolation of the functional $A_{N}^{{}}\xi =\sum\limits_{k=0}^{N}{a(k)\xi (k)}$ from unknown values of a stochastic sequence with stationary $n$th increments $\xi (m)$ based on observations of the sequence $\xi (k)$ at points $k=-1,-2,\ldots $, and observations of the sequence $\xi (k)+\eta (k)$ at points $k=N+1,N+2,\ldots $, where $\eta (k)$ is a stochastic sequence with stationary $n$ increments, uncorrelated with the sequence $\xi (k)$.
The problem of interpolation of a sequence with stationary $n$ increments is solved in the case where the spectral densities of the sequences $\xi (m)$ and $\eta (m)$ are known.
In the case where the spectral densities are unknown while some admissible sets of spectral densities are given, the minimax estimation method is applied.
For some given sets of admissible spectral densities, the least favorable densities and minimax spectral characteristics of the optimal estimate of the functional ${{A}_{N}}\xi $ are determined.

\noindent\hrulefill

\section{{Stationary increments. Spectral decomposition}}

\begin{ozn} \label{def1}
The stochastic $n$-th increment with step $\mu \in Z$ of a stochastic sequence $\{\xi (m),m\in Z\}$ is the function
\begin{equation}\label{01}
{{\xi }^{(n)}}(m,\mu )={{(1-{{B}_{\mu }})}^{n}}\xi (m)=\sum\limits_{l=0}^{n}{}\,{{(-1)}^{l}}C_{n}^{l}\xi (m-l\mu ),
\end{equation}
where ${{B}_{\mu }}$ is the shift operator by $\mu $ steps, such that ${{B}_{\mu }}\xi (m)=\xi (m-\mu )$, $m\in Z$.
\end{ozn}
For the stochastic $n$-th increment, the relations hold true
\begin{equation}\label{02}
{{\xi }^{(n)}}(m,-\mu )={{(-1)}^{n}}{{\xi }^{(n)}}(m+n\mu ,\mu ),
\end{equation}
\begin{equation}\label{04}
{{\xi }^{(n)}}(m,k\mu )=\sum\limits_{l=0}^{(k-1)n}{{{A}_{l}}{{\xi }^{(n)}}(m-l\mu ,\mu )}\,,\quad k\in N,
\end{equation}
where $\{{{A}_{l}},l=0,1,2,\ldots ,(k-1)n\}$ are the coefficients at the corresponding powers of ${{x}^{l}}$ of the polynomial expansion ${{(1+x+\ldots +{{x}^{k-1}})}^{n}}$ in powers of $x$.

\begin{ozn} \label{def2}
The stochastic $n$-th increment ${{\xi }^{(n)}}(m,\mu )$  with step $\mu \in Z$ of a stochastic sequence $\{\xi (m),m\in Z\}$  is called stationary (in the wide sense) if the mathematical expectations
\[\mathbf{E}{{\xi }^{(n)}}({{m}_{0}},\mu )={{c}^{(n)}}(\mu ),\quad \mathbf{E}{{\xi }^{(n)}}({{m}_{0}}+m,{{\mu }_{1}}){{\xi }^{(n)}}({{m}_{0}},{{\mu }_{2}})={{D}^{(n)}}(m,{{\mu }_{1}},{{\mu }_{2}})\]
exist for arbitrary integers ${{m}_{0}},$ $m,$ $\mu ,$ ${{\mu }_{1}},$ ${{\mu }_{2}}$, and do not depend on ${{m}_{0}}$.
The function ${{c}^{(n)}}(\mu )$ is called the mean value of the stationary $n$-th increment, and the function  ${{D}^{(n)}}(m,{{\mu }_{1}},{{\mu }_{2}})$ is called the structural function of the stationary $n$-th increment (or the structural function of the $n$-th order of the stochastic sequence $\{\xi (m),m\in Z\}$).\\
A stochastic sequence $\{\xi (m),m\in Z\}$ that determines the stationary $n$-th increment ${{\xi }^{(n)}}(m,\mu )$ by formula (1) is called a sequence with stationary $n$-th increments.
\end{ozn}

\begin{thm} \label{theorem01}
The mean value ${{c}^{(n)}}(\mu )$ and the structure function ${{D}^{(n)}}(m,{{\mu }_{1}},{{\mu }_{2}})$ of the stationary stochastic $n$-th increment ${{\xi }^{(n)}}(m,\mu )$ can be represented in the form
\begin{equation}\label{04}
{{c}^{(n)}}(\mu )=c{{\mu }^{n}},
\end{equation}
\begin{equation}\label{05}
{{D}^{(n)}}(m,{{\mu }_{1}},{{\mu }_{2}})=\int\limits_{-\pi }^{\pi }{{{e}^{i\lambda m}}{{(1-{{e}^{-i{{\mu }_{1}}\lambda }})}^{n}}{{(1-{{e}^{i{{\mu }_{2}}\lambda }})}^{n}}\frac{1}{{{\lambda }^{2n}}}dF(\lambda )}\,,
\end{equation}
where $c$ is a constant, the function $F(\lambda )$ is left continuous, nondecreasing, and bounded, $F(-\pi )=0$. The constant $c$ and the function $F(\lambda )$ are uniquely determined by the increment ${{\xi }^{(n)}}(m,\mu )$.\\
On the other hand, the function ${{c}^{(n)}}(\mu )$ of the form (4) with a constant $c$ and the function ${{D}^{(n)}}(m,{{\mu }_{1}},{{\mu }_{2}})$ of the form (5), where $F(\lambda )$ satisfies the above conditions, are the mean value and the structure function of some stationary $n$-th increment ${{\xi }^{(n)}}(m,\mu )$.
\end{thm}

Using the representation (5) of the structure function of the stationary $n$-th increment ${{\xi }^{(n)}}(m,\mu )$ and Karhunen's theorem [24], [25], we obtain the representation of the stationary $n$-th increment ${{\xi }^{(n)}}(m,\mu )$ in the following form
\begin{equation}\label{06}
{{\xi }^{(n)}}(m,\mu )=\int\limits_{-\pi }^{\pi }{{{e}^{im\lambda }}{{(1-{{e}^{-i\mu \lambda }})}^{n}}\frac{1}{{{(i\lambda )}^{n}}}dZ(\lambda )}\,,
\end{equation}
where $Z(\lambda )$ is an orthogonal stochastic measure on $[-\pi ,\pi )$, which is subordinate to the spectral measure generated by the function $F(\lambda )$:
\begin{equation}\label{07}
\mathbf{E}Z({{A}_{1}})\overline{Z({{A}_{2}})}=F({{A}_{1}}\cap {{A}_{2}})<\infty .
\end{equation}

We use the spectral decomposition (6) to find the optimal linear estimate of the unknown values of the stochastic sequence.

\section{{The interpolation problem}}

Let a stochastic sequence $\{\xi (m),m\in Z\}$ determine the stationary $n$-th increment ${{\xi }^{(n)}}(m,\mu )$ with an absolutely continuous spectral function $F(\lambda )$, which has the spectral density $f(\lambda )$, and the stochastic sequence $\{\eta (m),m\in Z\}$, which is uncorrelated with the sequence $\xi (m)$, determine the stationary $n$-th increment ${{\eta }^{(n)}}(m,\mu )$ with an absolutely continuous spectral function $G(\lambda )$, which has the spectral density $g(\lambda )$.

Suppose that we know the results of the observation of the sequence $\xi (k)$ at points $k=-1,-2,\ldots $ and the results of the observation of the sequence $\xi (k)+\eta (k)$ at points $k=N+1,N+2,\ldots $. Consider the problem of optimal in the mean square sense linear estimation of the functional
$${{A}_{N}}\xi =\sum\limits_{k=0}^{N}{a(k)\xi (k)}\,$$ from the unknown values of the sequence $\xi (m)$ based on the known observations.

From decomposition (1) we can obtain the following relation
\begin{equation}\label{08}
\xi (k)=\frac{1}{{{(1-{{B}_{\mu }})}^{n}}}{{\xi }^{(n)}}(k,\mu )=\sum\limits_{j=-\infty }^{k}{}\,{{d}_{\mu }}(k-j){{\xi }^{(n)}}(j,\mu ),
\end{equation}
where $\{{{d}_{\mu }}(k):k\ge 0\}$ are coefficients at ${{x}^{k}}$ in the decomposition
$$\sum\limits_{k=0}^{\infty }{{{d}_{\mu }}(k){{x}^{k}}}\,= {\left( \sum\limits_{k=0}^{\infty }{{{x}^{\mu j}}}\, \right)}^{n}.$$

From (1) and (8) the following relations follow
	\[\sum\limits_{k=0}^{N}{}\,a(k)\xi (k)=-\sum\limits_{i=-\mu n}^{-1}{}\,{{v}_{\mu }}(i)\xi (i)+\sum\limits_{i=0}^{N}{}\,\left( \sum\limits_{k=i}^{N}{}\,a(k){{d}_{\mu }}(k-i) \right){{\xi }^{(n)}}(i,\mu ),\]
	\[\sum\limits_{k=0}^{N}{}\,{{b}_{\mu }}(k){{\xi }^{(n)}}(k,\mu )=\sum\limits_{i=-\mu n}^{-1}{}\,\xi (i)\sum\limits_{l=\left[ -\tfrac{i}{\mu } \right]}^{n}{}\,{{(-1)}^{l}}C_{n}^{l}{{b}_{\mu }}(l\mu +i)+\sum\limits_{i=0}^{N}{}\,\xi (i)\sum\limits_{l=0}^{n}{}\,{{(-1)}^{l}}C_{n}^{l}{{b}_{\mu }}(l\mu +i),\]
where $[x{]}'$ is the smallest integer among the numbers greater than or equal to $x$, and ${{b}_{\mu }}(k)=0$ for $k>N$.
From the last two relations, we can obtain the following representation of the functional ${{A}_{N}}\xi $ in the form of a difference of functionals:
$${{A}_{N}}\xi ={{B}_{N}}\xi -{{V}_{N}}\xi,$$ where
	\[{{B}_{N}}\xi =\sum\limits_{k=0}^{N}{}\,{{b}_{\mu }}(k){{\xi }^{(n)}}(k,\mu ),\quad {{V}_{N}}\xi =\sum\limits_{k=-\mu n}^{-1}{}\,{{v}_{\mu }}(k)\xi (k),\]
\begin{equation}\label{09}
{{v}_{\mu }}(k)=\sum\limits_{l=\left[ -\tfrac{k}{\mu } \right]}^{n}{}\,{{(-1)}^{l}}C_{n}^{l}{{b}_{\mu }}(l\mu +k),\quad k=-1,-2,\ldots ,-\mu n,
\end{equation}
\begin{equation}\label{10}
{{b}_{\mu }}(k)=\sum\limits_{m=k}^{N}{}\,a(m){{d}_{\mu }}(m-k)={{(D_{N}^{\mu }{{a}^{(1)}})}_{k}},\,\,k=0,1,\ldots ,N.
\end{equation}
Here $D_{N}^{\mu }$ is the matrix of dimension $(N+1)\times (N+1)$ with elements $D_{k,j}^{\mu }={{d}_{\mu }}(j-k)$ if $0\le k\le j\le N$, and $D_{k,j}^{\mu }=0$ if $j<k$, $k,j=0,1,\ldots ,N$; vector ${{a}^{(1)}}=(a(0),a(1),a(2),\ldots ,a(N))$.

Let us denote by ${{\hat{A}}_{N}}\xi $ the optimal in the mean square sense linear estimate of the value of the functional ${{A}_{N}}\xi $ based on observations of random sequences $\xi (k)$ at points $k=-1,-2,\ldots $ and $\xi (k)+\eta (k)$ at points $k=N+1,N+2,\ldots $, and denote by ${{\hat{B}}_{N}}\xi $ the optimal in the mean square sense linear estimate of the value of the functional ${{B}_{N}}\xi $ based on observations of stochastic $n$-th increments ${{\xi }^{(n)}}(k,\mu )$ at points $k=-1,-2,\ldots $ and ${{\xi }^{(n)}}(k,\mu )+{{\eta }^{(n)}}(k,\mu )$ at points $k=N+1,N+2,\ldots $.

Let $\Delta (f,{{\hat{A}}_{N}})=\mathbf{E}|{{A}_{N}}\xi -{{\hat{A}}_{N}}\xi {{|}^{2}}$ denote the mean square error of the estimate ${{\hat{A}}_{N}}\xi $, and
let $\Delta (f,{{\hat{B}}_{N}})=\mathbf{E}|{{B}_{N}}\xi -{{\hat{B}}_{N}}\xi {{|}^{2}}$
denote the mean square error of the estimate ${{\hat{B}}_{N}}\xi $.
Since we know the values of the sequence $\xi (k)$ at points $k=-1,-2,\ldots ,-\mu n$, we can write the following equality:
\begin{equation}\label{11}
{{\hat{A}}_{N}}\xi ={{\hat{B}}_{N}}\xi -{{V}_{N}}\xi.
\end{equation}
Therefore, the following relations are valid:
	\[\Delta (f,{{\hat{A}}_{N}})=\mathbf{E}|{{A}_{N}}\xi -{{\hat{A}}_{N}}\xi {{|}^{2}}=\mathbf{E}|{{A}_{N}}\xi +{{V}_{N}}\xi -{{\hat{B}}_{N}}\xi {{|}^{2}}=\mathbf{E}|{{B}_{N}}\xi -{{\hat{B}}_{N}}\xi {{|}^{2}}=\Delta (f,{{\hat{B}}_{N}}).\]

To find the optimal linear estimate of the functional ${{B}_{N}}\xi $ in the mean-square sense, we will use the method of orthogonal projections in Hilbert space proposed by A.N. Kolmogorov [5].
Let us denote by ${{H}^{0-}}(\xi _{\mu }^{(n)})$ the closed linear subspace generated by $\{{{\xi }^{(n)}}(k,\mu ):k\le -1\}$ in the Hilbert space $H={{L}_{2}}(\Omega ,F,P)$ of second-order random variables, and by ${{H}^{N+}}(\xi _{-\mu }^{(n)}+\eta _{-\mu }^{(n)})$ we denote the closed linear subspace in the space $H={{L}_{2}}(\Omega ,F,P)$ generated by $\{{{\xi }^{(n)}}(k,-\mu )+{{\eta }^{(n)}}(k,-\mu ):k\ge N+1\}$.

Since
$${{e}^{i\lambda k}}{{(1-{{e}^{i\lambda \mu }})}^{n}}={{(-1)}^{n}}{{e}^{i\lambda (k+\mu n)}}{{(1-{{e}^{-i\lambda \mu }})}^{n}},$$
we have $${{\xi }^{(n)}+{{\eta }^{(n)}}(k,-\mu n)} )={{(-1)}^{n}}({{\xi }^{(n)}}(k+\mu n,\mu )+{{\eta }^{(n)}}(k+\mu n,\mu )).$$
 Therefore $${{H}^{N+}}(\xi _{-\mu }^{(n)}+\eta _{-\mu }^{(n)})={{H}^{(N+\mu n)+}}(\xi _{\mu }^{(n)}+\eta _{\mu }^{(n)}).$$
 We also define the subspaces $L_{2}^{0-}(f)$ and $L_{2}^{N+}(f+g)$ in the Hilbert spaces ${{L}_{2}}(f)$ and ${{L}_{2}}(f+g)$, which are generated by the functions
 $$\left\{{{e}^{i\lambda k}}{{(1-{{e}^{-i\lambda \mu }})}^{n}}\frac{1}{{{(i\lambda )}^{n}}}:k\le -1\right\}$$ and
 $$\left\{{{e}^{i\lambda k}}{{(1-{{e}^{-i\lambda \mu }})}^{n}}\frac{1}{{{(i\lambda )}^{n}}}:k\ge N+1\right\},$$ respectively.
We will search for a linear estimate of ${{\hat{B}}_{N}}\xi $ of the value of ${{B}_{N}}\xi $ in the form
\begin{equation}\label{12}
{{\hat{B}}_{N}}\xi =\int\limits_{-\pi }^{\pi }{h_{\mu }^{(1)}(\lambda )d{{Z}_{\xi _{\mu }^{(n)}}}(\lambda )}\,+
\int\limits_{-\pi }^{\pi }{h_{\mu }^{(2)}(\lambda )d{Z}_{{\xi _{\mu }^{(n)}}+\eta _{\mu }^{(n)}}}(\lambda ),
\end{equation}
where ${{h}_{\mu }}(\lambda )=(h_{\mu }^{(1)}(\lambda ),h_{\mu }^{(2)}(\lambda ))$ is the spectral characteristic of the estimate.
The optimal estimate ${{\hat{B}}_{N}}\xi $ is the projection of the element ${{B}_{N}}\xi $ onto the subspace $${{H}^{0-}}(\xi _{\mu }^{(n)})\oplus {{H}^{N+}}(\xi _{-\mu }^{(n)}+\eta _{-\mu }^{(n)})={{H}^{0-}}(\xi _{\mu }^{(n)})\oplus {{H}^{(N+\mu n)+}}(\xi _{\mu }^{(n)}+\eta _{\mu }^{(n)}).$$
Therefore, the estimate of $\widehat{B}\xi $ is determined by the following conditions:

1) ${{\hat{B}}_{N}}\xi \in \left({{H}^{0-}}(\xi _{\mu }^{(n)})\oplus {{H}^{(N+\mu n)+}}(\xi _{\mu }^{(n)}+\eta _{\mu }^{(n)})\right)$;

2) $\left({{B}_{N}}\xi -{{\hat{B}}_{N}}\xi\right)\perp \left({{H}^{0-}}(\xi _{\mu }^{(n)})\oplus {{H}^{(N+\mu n)+}}(\xi _{\mu }^{(n)}+\eta _{\mu }^{(n)})\right).$

\noindent From condition 2) it follows that for all $k\le -1$ the following relation holds
\[\mathbf{E}({{B}_{N}}\xi -{{\hat{B}}_{N}}\xi )\overline{{{\xi }^{(n)}}(k,\mu )}=\]
	\[=\frac{1}{2\pi }\int\limits_{-\pi }^{\pi }{\left( B_{N}^{\mu }({{e}^{i\lambda }}){{(1-{{e}^{-i\lambda \mu }})}^{n}}\frac{1}{{{(i\lambda )}^{n}}}-h_{\mu }^{(1)}(\lambda )-h_{\mu }^{(2)}(\lambda ) \right)\frac{{{e}^{-i\lambda k}}{{(1-{{e}^{i\lambda \mu }})}^{n}}}{{{(-i\lambda )}^{n}}}f(\lambda )d\lambda =0,}\,\]
and for all $k\ge N+\mu n+1$ the following relation holds
	\[\mathbf{E}({{B}_{N}}\xi -{{\hat{B}}_{N}}\xi )(\overline{{{\xi }^{(n)}}(k,\mu )+{{\eta }^{(n)}}(k,\mu )})=\]
	\[=\frac{1}{2\pi }\int\limits_{-\pi }^{\pi }{\left( B_{N}^{\mu }({{e}^{i\lambda }}){{(1-{{e}^{-i\lambda \mu }})}^{n}}\frac{1}{{{(i\lambda )}^{n}}}-h_{\mu }^{(1)}(\lambda )-h_{\mu }^{(2)}(\lambda ) \right)\frac{{{e}^{-i\lambda k}}{{(1-{{e}^{i\lambda \mu }})}^{n}}}{{{(-i\lambda )}^{n}}}f(\lambda )d\lambda -}\,\]
	\[-\frac{1}{2\pi }\int\limits_{-\pi }^{\pi }{\,h_{\mu }^{(2)}(\lambda ){{e}^{-i\lambda k}}{{(1-{{e}^{i\lambda \mu }})}^{n}}\frac{1}{{{(-i\lambda )}^{n}}}g(\lambda )d\lambda =}\]
	\[=\int\limits_{-\pi }^{\pi }{\,\left( B_{N}^{\mu }({{e}^{i\lambda }}){{(1-{{e}^{-i\lambda \mu }})}^{n}}\frac{f(\lambda )}{{{(i\lambda )}^{n}}}-h_{\mu }^{(1)}(\lambda )f(\lambda ) \right.}\left. -h_{\mu }^{(2)}(\lambda )(f(\lambda )+g(\lambda )) \right)\frac{{{e}^{-i\lambda k}}{{(1-{{e}^{i\lambda \mu }})}^{n}}}{{{(-i\lambda )}^{n}}}d\lambda =0.\]

Thus, the spectral characteristic ${{h}_{\mu }}(\lambda )=(h_{\mu }^{(1)}(\lambda ),h_{\mu }^{(2)}(\lambda ))$ of the estimate ${{\hat{B}}_{N}}\xi $ is given by the relations
	\[h_{\mu }^{(1)}(\lambda )+h_{\mu }^{(2)}(\lambda )=B_{N}^{\mu }({{e}^{i\lambda }}){{(1-{{e}^{-i\lambda \mu }})}^{n}}\frac{1}{{{(i\lambda )}^{n}}}-\frac{{{(-i\lambda )}^{n}}{{C}^{\mu }}({{e}^{i\lambda }})}{{{(1-{{e}^{i\lambda \mu }})}^{n}}f(\lambda )}\]
and
	\[h_{\mu }^{(2)}(\lambda )=B_{N}^{\mu }({{e}^{i\lambda }})\frac{{{(1-{{e}^{-i\lambda \mu }})}^{n}}}{{{(i\lambda )}^{n}}}\frac{f(\lambda )}{f(\lambda )+g(\lambda )}-\frac{h_{\mu }^{(1)}(\lambda )f(\lambda )}{f(\lambda )+g(\lambda )}-\frac{{{(-i\lambda )}^{n}}{{E}^{\mu }}({{e}^{i\lambda }})}{{{(1-{{e}^{i\lambda \mu }})}^{n}}(f(\lambda )+g(\lambda ))},\]
	\[B_{N}^{\mu }({{e}^{i\lambda }})=\sum\limits_{k=0}^{N}{}\,{{b}_{\mu }}(k){{e}^{i\lambda k}},\quad {{C}^{\mu }}({{e}^{i\lambda }})=\sum\limits_{k=0}^{\infty }{}\,{{c}_{\mu }}(k){{e}^{i\lambda k}},\quad {{E}^{\mu }}({{e}^{i\lambda }})=\sum\limits_{k=-\infty }^{N+\mu n}{}\,{{e}_{\mu }}(k){{e}^{i\lambda k}},\]
where ${{c}_{\mu }}(k)$, $k\ge 0$, and ${{e}_{\mu }}(k)$, $k\le N+\mu n$, are unknown coefficients.

From the last two equations we obtain the general form of the spectral characteristic
${{h}_{\mu }}(\lambda )=(h_{\mu }^{(1)}(\lambda ),h_{\mu }^{(2)}(\lambda ))$:
	\[h_{\mu }^{(1)}(\lambda )=B_{N}^{\mu }({{e}^{i\lambda }})\frac{{{(1-{{e}^{-i\lambda \mu }})}^{n}}}{{{(i\lambda )}^{n}}}+\frac{{{(-i\lambda )}^{n}}{{E}^{\mu }}({{e}^{i\lambda }})}{{{(1-{{e}^{i\lambda \mu }})}^{n}}g(\lambda )}-\frac{{{(-i\lambda )}^{n}}{{C}^{\mu }}({{e}^{i\lambda }})}{{{(1-{{e}^{i\lambda \mu }})}^{n}}}\left( \frac{1}{f(\lambda )}+\frac{1}{g(\lambda )} \right),\]
	\[h_{\mu }^{(2)}(\lambda )=({{C}^{\mu }}({{e}^{i\lambda }})-{{E}^{\mu }}({{e}^{i\lambda }}))\frac{{{(-i\lambda )}^{n}}}{{{(1-{{e}^{i\lambda \mu }})}^{n}}g(\lambda )}.\]

Let us find the equations that determine the coefficients ${{c}_{\mu }}(k)$, $k\ge 0$, and ${{e}_{\mu }}(k)$, $k\le N+\mu n$. From condition 1) it follows that the functions $h_{\mu }^{(1)}(\lambda )$ and $h_{\mu }^{(2)}(\lambda )$ have the form
	\[h_{\mu }^{(1)}(\lambda )={{h}^{(1)}}(\lambda ){{(1-{{e}^{-i\lambda \mu }})}^{n}}\frac{1}{{{(i\lambda )}^{n}}},\quad h_{\mu }^{(2)}(\lambda )={{h}^{(2)}}(\lambda ){{(1-{{e}^{-i\lambda \mu }})}^{n}}\frac{1}{{{(i\lambda )}^{n}}},\]
	\[{{h}^{(1)}}(\lambda )=\sum\limits_{k=-\infty }^{-1}{}\,s(k){{e}^{i\lambda k}},\quad {{h}^{(2)}}(\lambda )=\sum\limits_{k=N+\mu n+1}^{\infty }{}\,s(k){{e}^{i\lambda k}},\]
and satisfy the relation
	\[\int\limits_{-\pi }^{\pi }{\,|{{h}^{(1)}}(\lambda ){{|}^{2n}}|1-{{e}^{i\lambda \mu }}{{|}^{2}}f(\lambda ){{\lambda }^{2n}}d\lambda }<\infty ,\quad \int\limits_{-\pi }^{\pi }{\,|{{h}^{(2)}}(\lambda ){{|}^{2n}}|1-{{e}^{i\lambda \mu }}{{|}^{2}}f(\lambda )+g(\lambda ){{\lambda }^{2n}}d\lambda }<\infty, \]
	\[\frac{{{(i\lambda )}^{n}}h_{\mu }^{(1)}(\lambda )}{{{(1-{{e}^{-i\lambda \mu }})}^{n}}}\in L_{2}^{0-},\quad \frac{{{(i\lambda )}^{n}}h_{\mu }^{(2)}(\lambda )}{{{(1-{{e}^{-i\lambda \mu }})}^{n}}}\in L_{2}^{(N+\mu n)+},\]
\begin{equation}\label{13}
\int\limits_{-\pi }^{\pi }{\frac{{{(i\lambda )}^{n}}h_{\mu }^{(1)}(\lambda )}{{{(1-{{e}^{-i\lambda \mu }})}^{n}}}{{e}^{-i\lambda l}}d\lambda }\,=0,\quad l\ge 0.
\end{equation}
\begin{equation}\label{14}
\int\limits_{-\pi }^{\pi }{\frac{{{(i\lambda )}^{n}}h_{\mu }^{(2)}(\lambda )}{{{(1-{{e}^{-i\lambda \mu }})}^{n}}}{{e}^{-i\lambda l}}d\lambda =0}\,,\quad l\le N+\mu n.
\end{equation}
Let the following conditions be satisfied
\begin{equation}\label{15}
\int\limits_{-\pi }^{\pi } \frac{{\lambda }^{2n}}{|1-{{e}^{i\lambda \mu }}{{|}^{2n}} f(\lambda )}d\lambda \,<\infty,
\end{equation}	
\begin{equation}\label{16}
\int\limits_{-\pi }^{\pi } \frac{{\lambda }^{2n}}{|1-{{e}^{i\lambda \mu }}{{|}^{2n}} g(\lambda )}d\lambda \,<\infty.
\end{equation}
In this case we can determine the Fourier coefficients
\[{{f}_{\mu }}(k)=\frac{1}{2\pi }\int\limits_{-\pi }^{\pi }{\frac{{{\lambda }^{2n}}{{e}^{-i\lambda k}}}{|1-{{e}^{i\lambda \mu }}{{|}^{2n}}f(\lambda )}d\lambda }\,,\quad {{g}_{\mu }}(k)=\frac{1}{2\pi }\int\limits_{-\pi }^{\pi }{\frac{{{\lambda }^{2n}}{{e}^{-i\lambda k}}}{|1-{{e}^{i\lambda \mu }}{{|}^{2n}}g(\lambda )}d\lambda }\,,\quad k\in Z,\]
of the functions
$$\frac{{{\lambda }^{2n}}}{|1-{{e}^{i\lambda \mu }}{{|}^{2n}}f(\lambda )},\quad \frac{{{\lambda }^{2n}}}{|1-{{e}^{i\lambda \mu }}{{|}^{2n}}g(\lambda )}.$$
Suppose that the sequences $\xi (k)$ and $\eta (k)$ take real values. Then the properties ${{f}_{\mu }}(k)={{f}_{\mu }}(-k)$, ${{g}_{\mu }}(k)={{g}_{\mu }}(-k)$ hold true.
Using conditions (13) and (14) on the functions $h_{\mu }^{(1)}(\lambda )$ and $h_{\mu }^{(2)}(\lambda )$, we determine the unknown coefficients ${{c}_{\mu }}(k)$, $k\ge 0$, and ${{e}_{\mu }}(k)$, $k\le N+\mu n$. From condition (14) we obtain
	\[\int\limits_{-\pi }^{\pi }{\left( \sum\limits_{k=0}^{\infty }{}\,{{c}_{\mu }}(k){{e}^{i\lambda k}}-\sum\limits_{m=-\infty }^{N+\mu n}{}\,{{e}_{\mu }}(m){{e}^{i\lambda m}} \right)\frac{{{\lambda }^{2n}}{{e}^{-i\lambda l}}}{|1-{{e}^{-i\lambda \mu }}{{|}^{2n}}g(\lambda )}d\lambda =0}\,,\quad l\le N+\mu n,\]
	\[\sum\limits_{k=0}^{\infty }{}\,{{c}_{\mu }}(k)g(k-l)=\sum\limits_{m=-\infty }^{N+\mu n}{}\,{{e}_{\mu }}(m)g(l-m),\quad l\le N+\mu n.\]
The last system can be written as
	\[\mathbf{G}_{\mu }^{c}{{\mathbf{c}}_{\mu }}=\mathbf{G}_{\mu }^{e}{{\mathbf{e}}_{\mu }},\]
where ${{\mathbf{c}}_{\mu }}=({{c}_{\mu }}(0),{{c}_{\mu }}(1),{{c}_{\mu }}(2),\ldots )$, ${{\mathbf{e}}_{\mu }}=({{e}_{\mu }}(N+\mu n),{{e}_{\mu }}(N+\mu n-1),{{e}_{\mu }}(N+\mu n-2),\ldots )$, $\mathbf{G}_{\mu }^{c}$, $\mathbf{G}_{\mu }^{e}$ are linear operators in the space ${{\ell}_{2}}$, defined by matrices with elements
$${{(\mathbf{G}_{\mu }^{c})}_{l,k}}={{g}_{\mu }}(N+\mu n-l-k), \quad {{(\mathbf{G}_{\mu }^{e})}_{l,k}}={{g}_{\mu }}(k-l),\,\, l,k\ge 0,$$ respectively.
From condition (13) we obtain the following relations:
	\[\int\limits_{-\pi }^{\pi }{\left( {{B}_{\mu }}({{e}^{i\lambda }})-\sum\limits_{k=0}^{\infty }{}\,{{c}_{\mu }}(k){{e}^{i\lambda k}}\frac{{{\lambda }^{2n}}}{|1-{{e}^{-i\lambda \mu }}{{|}^{2n}}g(\lambda )}\left( \frac{1}{f(\lambda )}+\frac{1}{g(\lambda )} \right)+ \right.}\,\]
	\[\left. +\sum\limits_{m=-\infty }^{N+\mu n}{}\,{{e}_{\mu }}(m){{e}^{i\lambda m}}\frac{{{\lambda }^{2n}}}{|1-{{e}^{-i\lambda \mu }}{{|}^{2n}}g(\lambda )} \right){{e}^{-i\lambda l}}d\lambda =0,\quad l\ge 0,\]
	\[{{b}_{\mu }}(l)+\sum\limits_{k=0}^{\infty }{}\,{{c}_{\mu }}(k){{g}_{\mu }}(l-k)=\sum\limits_{m=-\infty }^{N+\mu n}{}\,{{e}_{\mu }}(m)({{g}_{\mu }}(l-m)+{{f}_{\mu }}(l-m)),\quad l\ge 0.\]
Using the properties ${{f}_{\mu }}(k)={{f}_{\mu }}(-k)$ and ${{g}_{\mu }}(k)={{g}_{\mu }}(-k)$, from the last system of equations we can obtain the following relation, which holds for any $l\ge 0$:
	\[{{b}_{\mu }}(l)+\sum\limits_{k=0}^{\infty }{}\,{{c}_{\mu }}(k){{g}_{\mu }}(k-l)=\sum\limits_{m=0}^{\infty }{}\,{{e}_{\mu }}(N+\mu n-m)({{g}_{\mu }}(N+\mu n-l-m)+{{f}_{\mu }}(N+\mu n-l-m)).\]
The resulting ratios are written in the form
	\[{{\mathbf{b}}_{\mu ,N}}+\mathbf{G}_{\mu }^{e}{{\mathbf{c}}_{\mu }}=\mathbf{G}_{\mu }^{c}{{\mathbf{e}}_{\mu }}+\mathbf{F}_{\mu }^{c}{{\mathbf{e}}_{\mu }},\]
where
${{\mathbf{b}}_{\mu ,N}}=({{b}_{\mu }}(0),{{b}_{\mu }}(1),\ldots ,{{b}_{\mu }}(N),0,\ldots )$, and $\mathbf{F}_{\mu }^{c}$
is the linear operator in the space ${{\ell}_{2}}$, defined by the matrix with elements $${{(\mathbf{F}_{\mu }^{c})}_{l,k}}={{f}_{\mu }}(N+\mu n-l-k),\quad l,k\ge 0.$$

Thus, the coefficients ${{c}_{\mu }}(k)$, $k\ge 0$, and ${{e}_{\mu }}(k)$, $k\le N+\mu n$ are determined by the following systems:
	\[{{\mathbf{b}}_{\mu ,N}}+\mathbf{G}_{\mu }^{e}{{\mathbf{c}}_{\mu }}=\mathbf{G}_{\mu }^{c}{{\mathbf{e}}_{\mu }}+\mathbf{F}_{\mu }^{c}{{\mathbf{e}}_{\mu }},\quad \mathbf{G}_{\mu }^{c}{{\mathbf{c}}_{\mu }}=\mathbf{G}_{\mu }^{e}{{\mathbf{e}}_{\mu }}.\]
Suppose that the matrices
$$\mathbf{G}_{\mu }^{e},\,\,{{\mathbf{F}}_{\mu }}=\mathbf{G}_{\mu }^{c}{{(\mathbf{G}_{\mu }^{e})}^{-1}}\mathbf{G}_{\mu }^{c}+\mathbf{F}_{\mu }^{c}{{(\mathbf{G}_{\mu }^{e})}^{-1}}\mathbf{G}_{\mu }^{c}-\mathbf{G}_{\mu }^{e}$$
are invertible. Then the unknown coefficients are determined as follows:
	\[{{\mathbf{c}}_{\mu }}={{({{\mathbf{F}}_{\mu }})}^{-1}}{{\mathbf{b}}_{\mu ,N}},\quad {{\mathbf{e}}_{\mu }}={{({{\mathbf{G}}_{\mu }})}^{-1}}{{\mathbf{b}}_{\mu ,N}},\]
where
	\[{{\mathbf{F}}_{\mu }}=\mathbf{G}_{\mu }^{c}{{(\mathbf{G}_{\mu }^{e})}^{-1}}\mathbf{G}_{\mu }^{c}+\mathbf{F}_{\mu }^{c}{{(\mathbf{G}_{\mu }^{e})}^{-1}}\mathbf{G}_{\mu }^{c}-\mathbf{G}_{\mu }^{e},\quad {{\mathbf{G}}_{\mu }}={{(\mathbf{G}_{\mu }^{e})}^{-1}}\mathbf{G}_{\mu }^{c}{{({{\mathbf{F}}_{\mu }})}^{-1}},\]
${{e}_{\mu }}(k)={{({{({{\mathbf{G}}_{\mu }})}^{-1}}{{\mathbf{b}}_{\mu ,N}})}_{N+\mu n-k}}$ is $(N+\mu n-k)$-th element of the vector ${{({{\mathbf{G}}_{\mu }})}^{-1}}{{\mathbf{b}}_{\mu ,N}}$, $k\le N+\mu n$; ${{c}_{\mu }}(k)={{({{({{\mathbf{F}}_{\mu }})}^{-1}}{{\mathbf{b}}_{\mu ,N}})}_{k}}$ is $k$-th element of the vector ${{({{\mathbf{F}}_{\mu }})}^{-1}}{{\mathbf{b}}_{\mu ,N}}$, $k\ge 0$.

Thus, the spectral characteristic ${{h}_{\mu }}(\lambda )=(h_{\mu }^{(1)}(\lambda ),h_{\mu }^{(2)}(\lambda ))$ of the optimal estimate of the ${{\hat{B}}_{N}}\xi $ functional is calculated by the formulas
	\[h_{\mu }^{(1)}(\lambda )=B_{N}^{\mu }({{e}^{i\lambda }})\frac{{{(1-{{e}^{-i\lambda \mu }})}^{n}}}{{{(i\lambda )}^{n}}}+\frac{{{(-i\lambda )}^{n}}\sum\limits_{k=-\infty }^{N+\mu n}{}\,{{({{({{\mathbf{G}}_{\mu }})}^{-1}}{{\mathbf{b}}_{\mu ,N}})}_{N+\mu n-k}}{{e}^{i\lambda k}}}{{{(1-{{e}^{i\lambda \mu }})}^{n}}g(\lambda )}-\]
\begin{equation}\label{17}
-\frac{{{(-i\lambda )}^{n}}\sum\limits_{k=0}^{\infty }{}\,{{({{({{\mathbf{F}}_{\mu }})}^{-1}}{{\mathbf{b}}_{\mu ,N}})}_{k}}{{e}^{i\lambda k}}}{{{(1-{{e}^{i\lambda \mu }})}^{n}}}\left( \frac{1}{f(\lambda )}+\frac{1}{g(\lambda )} \right),
\end{equation}
\begin{equation}\label{18}
h_{\mu }^{(2)}(\lambda )=\left( \sum\limits_{k=0}^{\infty }{}\,{{({{({{\mathbf{F}}_{\mu }})}^{-1}}{{\mathbf{b}}_{\mu ,N}})}_{k}}{{e}^{i\lambda k}}-\sum\limits_{k=-\infty }^{N+\mu n}{}\,{{({{({{\mathbf{G}}_{\mu }})}^{-1}}{{\mathbf{b}}_{\mu ,N}})}_{N+\mu n-k}}{{e}^{i\lambda k}} \right)\frac{{{(-i\lambda )}^{n}}}{{{(1-{{e}^{i\lambda \mu }})}^{n}}g(\lambda )}.
\end{equation}
The mean-square error of the estimate is calculated by the formula
	\[\Delta (f,g;{{\hat{B}}_{N}})=\int\limits_{-\pi }^{\pi }{\frac{{{\lambda }^{2n}} {\left| \sum\limits_{k=0}^{\infty }{}{{({{({{\mathbf{F}}_{\mu }})}^{-1}}{{\mathbf{b}}_{\mu ,N}})}_{k}}{{e}^{i\lambda k}} \right|}^{2}}{|1-{{e}^{i\lambda \mu }}{{|}^{2n}}f(\lambda )}d\lambda }\,+\]
\begin{equation}\label{19}
+\int\limits_{-\pi }^{\pi }{\frac{{{\lambda }^{2n}} {\left| \sum\limits_{k=0}^{\infty }{}{{({{({{\mathbf{F}}_{\mu }})}^{-1}}{{\mathbf{b}}_{\mu ,N}})}_{k}}{{e}^{i\lambda k}}-\sum\limits_{k=-\infty }^{N+\mu n}{}\,{{({{({{\mathbf{G}}_{\mu }})}^{-1}}{{\mathbf{b}}_{\mu ,N}})}_{N+\mu n-k}}{{e}^{i\lambda k}} \right|}^{2}}{|1-{{e}^{i\lambda \mu }}{{|}^{2n}}g(\lambda )}d\lambda }\,.
\end{equation}

Therefore, the following theorem is valid.

\begin{thm} \label{theorem02}
 Let the uncorrelated stochastic sequences $\{\xi (m),m\in Z\}$ and $\{\eta (m),m\in Z\}$ determine the stationary $n$-th increments ${{\xi }^{(n)}}(m,\mu )$ and ${{\eta }^{(n)}}(m,\mu )$ with absolutely continuous spectral functions $F(\lambda )$ and $G(\lambda )$, which have spectral densities $f(\lambda )$ and $g(\lambda )$ satisfying conditions (15) and (16). Then the optimal linear estimate ${{\hat{B}}_{N}}\xi $ of the value of the functional ${{B}_{N}}\xi $ from unknown increments ${{\xi }^{(n)}}(m,\mu )$, $m\in \{0,1,2,\ldots ,N\}$, $\mu >0$, based on observations of the sequence $\xi (m)$ at points $m=-1,-2,\ldots $ and observations of the sequence $\xi (m)+\eta (m)$ at points $m=N+1,N+2,\ldots $ is determined by formula (12).
 The spectral characteristic ${{h}_{\mu }}(\lambda )=(h_{\mu }^{(1)}(\lambda ),h_{\mu }^{(2)}(\lambda ))$ of the optimal estimate ${{\hat{B}}_{N}}\xi $ is calculated by formulas (17) and (18).
 The value of the mean square error $\Delta (f,{{\hat{B}}_{N}})$ of the estimate ${{\hat{B}}_{N}}\xi $ is calculated by formula (19).
\end{thm}

We use this Theorem \ref{theorem02} to construct an estimate of the unknown value of the increment ${{\xi }^{(n)}}(m,\mu )$, $m\in \{0,1,\ldots ,N\}$, based on observations of the sequence $\xi (m)$ at points $m=-1,-2,\ldots $ and the sequence $\xi (m)+\eta (m)$ at points $m=N+1,N+2,\ldots $ . We substitute the vector ${{b}_{\mu ,m}}$ into formulas (17) and (18), in which the $m$-th place is $1$, and all other elements are zero. Then the formulas for calculating the spectral characteristic ${{\varphi }_{m}}(\lambda ,\mu )=(\varphi _{m}^{(1)}(\lambda ,\mu ),\varphi _{m}^{(2)}(\lambda ,\mu ))$ of the estimate
\begin{equation}\label{20}
{{\hat{\xi }}}^{(n)}(m,\mu )=\int\limits_{-\pi }^{\pi }{\varphi _{m}^{(1)}(\lambda ,\mu )d{{Z}_{\xi _{\mu }^{(n)}}}(\lambda )}\,+\int\limits_{-\pi }^{\pi }{\varphi _{m}^{(2)}(\lambda ,\mu )d{{Z}_{\xi _{\mu }^{(n)}+\eta _{\mu }^{(n)}}}(\lambda )}\,
\end{equation}
will take the form
\begin{multline}\label{21}
\varphi _{m}^{(1)}(\lambda ,\mu )={{e}^{i\lambda m}}\frac{{{(1-{{e}^{-i\lambda \mu }})}^{n}}}{{{(i\lambda )}^{n}}}+\\
+\frac{{{e}^{i\lambda (N+\mu n)}}\sum\limits_{k=0}^{\infty }{}\,{{({{\mathbf{W}}_{\mu }})}_{k,m}}{{e}^{-i\lambda k}}}{{{(-i\lambda )}^{-n}}{{(1-{{e}^{i\lambda \mu }})}^{n}}g(\lambda )}-\frac{\sum\limits_{k=0}^{\infty }{}{{({{\mathbf{V}}_{\mu }})}_{k,m}}{{e}^{i\lambda k}}}{{{(-i\lambda )}^{-n}}{{(1-{{e}^{i\lambda \mu }})}^{n}}}\left( \frac{1}{f(\lambda )}+\frac{1}{g(\lambda )} \right),
\end{multline}
\begin{equation}\label{22}
\varphi _{m}^{(2)}(\lambda ,\mu )=\left( \sum\limits_{k=0}^{\infty }{}\,{{({{\mathbf{V}}_{\mu }})}_{k,m}}{{e}^{i\lambda k}}-{{e}^{i\lambda (N+\mu n)}}\sum\limits_{k=0}^{\infty
}{}\,{{({{\mathbf{W}}_{\mu }})}_{k,m}}{{e}^{-i\lambda k}} \right)\frac{{{(-i\lambda )}^{n}}}{{{(1-{{e}^{i\lambda \mu }})}^{n}}g(\lambda )},
\end{equation}
where ${{\mathbf{V}}_{\mu }}={{({{\mathbf{F}}_{\mu }})}^{-1}}$, ${{\mathbf{W}}_{\mu }}={{({{\mathbf{G}}_{\mu }})}^{-1}}$.
\\
The value of the mean square error is calculated by formula
\begin{multline}\label{23}
\Delta (f,{{\hat{\xi }}^{(n)}}(m,\mu ))=\int\limits_{-\pi }^{\pi }{\frac{ {\left| \sum\limits_{k=0}^{\infty }{}{{({{\mathbf{V}}_{\mu }})}_{k,m}}{{e}^{i\lambda k}} \right|}^{2}}{{{\lambda }^{-2n}}|1-{{e}^{i\lambda \mu }}{{|}^{2n}}f(\lambda )}d\lambda }\,+\\
+\int\limits_{-\pi }^{\pi }{\frac{ {\left| \sum\limits_{k=0}^{\infty }{}{{({{\mathbf{V}}_{\mu }})}_{k,m}}{{e}^{i\lambda k}}-{{e}^{i\lambda (N+\mu n)}}\sum\limits_{k=0}^{\infty }{}\,{{({{\mathbf{W}}_{\mu }})}_{k,m}}{{e}^{-i\lambda k}} \right|}^{2}}{{{\lambda }^{-2n}}|1-{{e}^{i\lambda \mu }}{{|}^{2n}}g(\lambda )}d\lambda }\,.
\end{multline}

Thus, the following corollary holds true.

\begin{nas} \label{corr1}
 The optimal linear estimate ${{\hat{\xi }}^{(n)}}(m,\mu )$ of the unknown value of the increment ${{\xi }^{(n)}}(m,\mu )$, $m=0,1,2,\ldots ,N$, $\mu >0$, based on the observations of the sequence $\xi (m)$ at points $m=-1,-2,\ldots $ and observations of the sequence $\xi (m)+\eta (m)$ at points $m=N+1,N+2,\ldots $ is determined by formula (20). The spectral characteristic $\varphi (\lambda ,\mu )$ of the optimal estimate ${{\hat{\xi }}^{(n)}}(m,\mu )$ is calculated by formulas (21) and (22). The value of the mean square error $\Delta (f,{{\hat{\xi }}^{(n)}}(m,\mu ))$ is calculated by formula (23).
\end{nas}

Using formula (11) and this Theorem \ref{theorem02}, we construct an estimate
\begin{equation}\label{24}
{{\hat{A}}_{N}}\xi =-\sum\limits_{k=-\mu n}^{-1}{}\,{{v}_{\mu }}(k)\xi (k)+\int\limits_{-\pi }^{\pi }{h_{\mu }^{(a,1)}(\lambda )d{{Z}_{\xi _{\mu }^{(n)}}}(\lambda )}\,+\int\limits_{-\pi }^{\pi }{h_{\mu }^{(a,2)}(\lambda )d{{Z}_{\xi _{\mu }^{(n)}+\eta _{\mu }^{(n)}}}(\lambda )}\,
\end{equation}
of the functional ${{A}_{N}}\xi $. Substitute the vector ${{\mathbf{b}}_{\mu ,N}}$ with the coefficients defined in (10) into relations (17), (18) and (19). Denote by ${{[D_{N}^{\mu }{{a}^{(1)}}]}_{\infty }}$ the embedding of the $(N+1)$-dimensional vector $D_{N}^{\mu }{{a}^{(1)}}$ into the space ${{\ell}_{2}}$ such that the first $N+1$ coordinates of the element ${{[D_{N}^{\mu }{{a}^{(1)}}]}_{\infty }}$ are equal to the coordinates of the vector $D_{N}^{\mu }{{a}^{(1)}}$, and all other coordinates are equal to $0$. We obtain the following formulas for determination the spectral characteristic $h_{\mu }^{(a)}(\lambda )=(h_{\mu }^{(a,1)}(\lambda ),h_{\mu }^{(a,2)}(\lambda ))$ and the mean square error $\Delta (f,{{\hat{A}}_{N}})$ of the estimate ${{\hat{A}}_{N}}\xi $:
\begin{multline}\label{25}
h_{\mu }^{(a,1)}(\lambda )=A_{N}^{\mu }({{e}^{i\lambda }})\frac{{{(1-{{e}^{-i\lambda \mu }})}^{n}}}{{{(i\lambda )}^{n}}}+\frac{{{(-i\lambda )}^{n}}\sum\limits_{k=-\infty }^{N+\mu n}{}\,{{({{({{\mathbf{G}}_{\mu }})}^{-1}}{{[D_{N}^{\mu }{{a}^{(1)}}]}_{\infty }})}_{N+\mu n-k}}{{e}^{i\lambda k}}}{{{(1-{{e}^{i\lambda \mu }})}^{n}}g(\lambda )}-\\
	-\frac{{{(-i\lambda )}^{n}}\sum\limits_{k=0}^{\infty }{}\,{{({{({{\mathbf{F}}_{\mu }})}^{-1}}{{[D_{N}^{\mu }{{a}^{(1)}}]}_{\infty }})}_{k}}{{e}^{i\lambda k}}}{{{(1-{{e}^{i\lambda \mu }})}^{n}}}\left( \frac{1}{f(\lambda )}+\frac{1}{g(\lambda )} \right),
\end{multline}

\begin{multline}\label{26}
	h_{\mu }^{(a,2)}(\lambda )=\frac{{{(-i\lambda )}^{n}}\sum\limits_{k=0}^{\infty }{}\,{{({{({{\mathbf{F}}_{\mu }})}^{-1}}{{[D_{N}^{\mu }{{a}^{(1)}}]}_{\infty }})}_{k}}{{e}^{i\lambda k}}}{{{(1-{{e}^{i\lambda \mu }})}^{n}}g(\lambda )}-\\
-\frac{{{(-i\lambda )}^{n}}\sum\limits_{k=-\infty }^{N+\mu n}{}\,{{({{({{\mathbf{G}}_{\mu }})}^{-1}}{{[D_{N}^{\mu }{{a}^{(1)}}]}_{\infty }})}_{N+\mu n-k}}{{e}^{i\lambda k}}}{{{(1-{{e}^{i\lambda \mu }})}^{n}}g(\lambda )},
\end{multline}

where $A_{N}^{\mu }({{e}^{i\lambda }})=\sum\limits_{k=0}^{N}{\,{{(D_{N}^{\mu }{{a}^{(1)}})}_{k}}{{e}^{i\lambda k}}}$,

\begin{multline}\label{27}
	\Delta (f,g;{{\hat{A}}_{N}})=\int\limits_{-\pi }^{\pi }{\frac{ {\left| \sum\limits_{k=0}^{\infty }{}{{({{({{\mathbf{F}}_{\mu }})}^{-1}}{{[D_{N}^{\mu }{{a}^{(1)}}]}_{\infty }})}_{k}}{{e}^{i\lambda k}} \right|}^{2}}{{{\lambda }^{-2n}}|1-{{e}^{i\lambda \mu }}{{|}^{2n}}f(\lambda )}d\lambda }\,+
\\
	+\int\limits_{-\pi }^{\pi }{\,\frac{ {\left| \sum\limits_{k=0}^{\infty }{}{{({{({{\mathbf{F}}_{\mu }})}^{-1}}{{[D_{N}^{\mu }{{a}^{(1)}}]}_{\infty }})}_{k}}{{e}^{i\lambda k}}-\sum\limits_{k=-\infty }^{N+\mu n}{}\,{{({{({{\mathbf{G}}_{\mu }})}^{-1}}{{[D_{N}^{\mu }{{a}^{(1)}}]}_{\infty }})}_{N+\mu n-k}}{{e}^{i\lambda k}} \right|}^{2}}{{{\lambda }^{-2n}}|1-{{e}^{i\lambda \mu }}{{|}^{2n}}g(\lambda )}d\lambda .}
\end{multline}

Therefore, the following theorem holds true.

\begin{thm} \label{theorem03}
Let the uncorrelated stochastic sequences $\{\xi (m),m\in Z\}$ and $\{\eta (m),m\in Z\}$ determine the stationary $n$-th increments ${{\xi }^{(n)}}(m,\mu )$ and ${{\eta }^{(n)}}(m,\mu )$ with absolutely continuous spectral functions $F(\lambda )$ and $G(\lambda )$, which have spectral densities $f(\lambda )$ and $g(\lambda )$ satisfying conditions (15) and (16). Then the optimal linear estimate ${{\hat{A}}_{N}}\xi $ of the value of the functional ${{A}_{N}}\xi $ from the unknown values $\xi (m)$, $m\in \{0,1,2,\ldots ,N\}$, based on the observations of the sequence $\xi (m)$ at points $m=-1,-2,\ldots $ and observations of the sequence $\xi (m)+\eta (m)$ at points $m=N+1,N+2,\ldots $ is determined by formula (24).
The spectral characteristic $h_{\mu }^{(a)}(\lambda )=(h_{\mu }^{(a,1)}(\lambda ),h_{\mu }^{(a,2)}(\lambda ))$ of the optimal estimate ${{\hat{A}}_{N}}\xi $ is calculated by formulas (25) and (26).
The value of the mean square error $\Delta (f,{{\hat{A}}_{N}})$ of the estimate ${{\hat{A}}_{N}}\xi $ is calculated by formula (27).
\end{thm}

\section{{The relationship between of the problem of interpolation from observations with noise and the problem of filtering}}

When constructing the estimate ${{\hat{A}}_{N}}\xi $ of the functional ${{A}_{N}}\xi $, we find a projection of the functional ${{B}_{N}}\xi $ onto the space ${{H}^{0-}}(\xi _{\mu }^{(n)})\oplus {{H}^{(N+\mu n)+}}(\xi _{\mu }^{(n)}+\eta _{\mu }^{(n)})$.
However, the problem can be generalized by studying the problem of estimating the functional

	\[{{B}_{N+\mu n}}\xi =\sum\limits_{k=0}^{N+\mu n}{}\,{{b}_{\mu }}(k){{\xi }^{(n)}}(k,\mu ).\]
Using relation (11), we can construct an estimate of the functional
	\[{{A}_{N+\mu n}}\xi =\sum\limits_{k=0}^{N+\mu n}{}\,a(k)\xi (k),\]
replacing $N$ with $N+\mu n$. If we put $a(k)=0$ for $k=0,1,\ldots ,N$, then we obtain the problem of filtering the functional
	\[A_{N+1}^{N+\mu n}\xi =\sum\limits_{k=N+1}^{N+\mu n}{a(k)\xi (k)}\]
from unknown values $\xi (k)$, $k=N+1,N+2,\ldots ,N+\mu n,$ of the sequence $\xi (k)$ based on observations of the sequence $\xi (k)$ at points $k=-1,-2,\ldots $ and observations of the sequence $\xi (k)+\eta (k)$ at points $k=N+1,N+2,\ldots $. The estimate $\hat{A}_{N+1}^{N+\mu n}\xi $ of the functional $A_{N+1}^{N+\mu n}\xi $ can be represented in the form $$\hat{A}_{N+1}^{N+\mu n}\xi ={{\hat{B}}_{N+\mu n}}\xi -{{V}_{N+\mu n}}\xi,$$ where
\begin{equation}\label{28}
{{V}_{N+\mu n}}\xi =\sum\limits_{k=-\mu n}^{-1}{}\,{{\tilde{v}}_{\mu }}(k)\xi (k),\quad {{\tilde{v}}_{\mu }}(k)=\sum\limits_{l=\left[ -\tfrac{k}{\mu } \right]}^{n}{}\,{{(-1)}^{l}}C_{n}^{l}{{\tilde{b}}_{\mu }}(l\mu +k),\quad k=-1,-2,\ldots ,-\mu n,
\end{equation}
\begin{equation}\label{29}
{{\tilde{b}}_{\mu }}(k)=\sum\limits_{m=max\{k,N+1\}}^{N+\mu n}{}\,a(m){{d}_{\mu }}(m-k)={{(D_{N+\mu n}^{\mu }{{a}^{(2)}})}_{k}},\quad k=0,1,\ldots ,N+\mu n,
\end{equation}
${{\tilde{b}}_{\mu }}(k)=0$, $k>N+\mu n$. The matrix $D_{N+\mu n}^{\mu }$ in the latter relation has dimension $(N+\mu n+1)\times (N+\mu n+1)$ and is defined by the elements $D_{k,j}^{\mu }={{d}_{\mu }}(j-k)$, if $0\le k\le j\le N+\mu n$, and $D_{k,j}^{\mu }=0$, if $j<k$, $k,j=0,1,\ldots ,N+\mu n$; the vector ${{a}^{(2)}}=(0,0,\ldots ,0,a(N+1),a(N+2),\ldots ,a(N+\mu n))$ has dimension $(N+\mu n+1)$.

Therefore, using the proved Theorem \ref{theorem02} to construct the estimate ${{\hat{B}}_{N+\mu n}}\xi $, the estimate $\hat{A}_{N+1}^{N+\mu n}\xi $ of the functional $A_{N+1}^{N+\mu n}\xi $ can be represented as
\begin{equation}\label{30}
A_{N+1}^{N+\mu n}\xi =-\sum\limits_{k=-\mu n}^{-1}{}\,\,{{\tilde{v}}_{\mu }}(k)\xi (k)+\int\limits_{-\pi }^{\pi }{\tilde{h}_{\mu }^{(a,1)}(\lambda )d{{Z}_{\xi _{\mu }^{(n)}}}(\lambda )}\,+\int\limits_{-\pi }^{\pi }{\tilde{h}_{\mu }^{(a,2)}(\lambda )d{{Z}_{\xi _{\mu }^{(n)}+\eta _{\mu }^{(n)}}}(\lambda ),}\,
\end{equation}
where the spectral characteristic $\tilde{h}_{\mu }^{(a)}(\lambda )=(\tilde{h}_{\mu }^{(a,1)}(\lambda ),\tilde{h}_{\mu }^{(a,2)}(\lambda ))$ is calculated by the formulas
\begin{multline}\label{31}
\tilde{h}_{\mu }^{(a,1)}(\lambda )=A_{N+1,\mu }^{N+\mu n}({{e}^{i\lambda }})\frac{{{(1-{{e}^{-i\lambda \mu }})}^{n}}}{{{(i\lambda )}^{n}}}+\frac{\sum\limits_{k=-\infty }^{N+\mu n}{}\,{{({{({{\mathbf{G}}_{\mu }})}^{-1}}{{[D_{N+\mu n}^{\mu }{{a}^{(2)}}]}_{\infty }})}_{N+\mu n-k}}{{e}^{i\lambda k}}}{{{(-i\lambda )}^{-n}}{{(1-{{e}^{i\lambda \mu }})}^{n}}g(\lambda )}-
\\
-\frac{\sum\limits_{k=0}^{\infty }{}\,{{({{({{\mathbf{F}}_{\mu }})}^{-1}}{{[D_{N+\mu n}^{\mu }{{a}^{(2)}}]}_{\infty }})}_{k}}{{e}^{i\lambda k}}}{{{(-i\lambda )}^{-n}}{{(1-{{e}^{i\lambda \mu }})}^{n}}}\left( \frac{1}{f(\lambda )}+\frac{1}{g(\lambda )} \right),
\end{multline}
\begin{multline}\label{32}
\tilde{h}_{\mu }^{(a,2)}(\lambda )=\frac{\sum\limits_{k=0}^{\infty }{}\,{{({{({{\mathbf{F}}_{\mu }})}^{-1}}{{[D_{N+\mu n}^{\mu }{{a}^{(2)}}]}_{\infty }})}_{k}}{{e}^{i\lambda k}}}{{{(-i\lambda )}^{-n}}{{(1-{{e}^{i\lambda \mu }})}^{n}}g(\lambda )}-
\\
-\frac{\sum\limits_{k=-\infty }^{N+\mu n}{}\,{{({{({{\mathbf{G}}_{\mu }})}^{-1}}{{[D_{N+\mu n}^{\mu }{{a}^{(1)}}]}_{\infty }})}_{N+\mu n-k}}{{e}^{i\lambda k}}}{{{(-i\lambda )}^{-n}}{{(1-{{e}^{i\lambda \mu }})}^{n}}g(\lambda )},
\end{multline}

$A_{N+1,\mu }^{N+\mu n}({{e}^{i\lambda }})=\sum\limits_{k=N+1}^{N+\mu n}{{{(D_{N+\mu n}^{\mu }{{a}^{(2)}})}_{k}}{{e}^{i\lambda k}}}.$\\
The value of the mean square error $\Delta (f,\hat{A}_{N+1}^{N+\mu n})$ of the estimate $\hat{A}_{N+1}^{N+\mu n}\xi $  is calculated by the formula
\begin{multline}\label{33}
	\Delta (f,g;\hat{A}_{N+1}^{N+\mu n})=\int\limits_{-\pi }^{\pi }{\frac{ {\left| \sum\limits_{k=0}^{\infty }{}{{({{({{\mathbf{F}}_{\mu }})}^{-1}}{{[D_{N+\mu n}^{\mu }{{a}^{(2)}}]}_{\infty }})}_{k}}{{e}^{i\lambda k}} \right|}^{2}}{{{\lambda }^{-2n}}|1-{{e}^{i\lambda \mu }}{{|}^{2n}}f(\lambda )}d\lambda }\,+
\\
+\int\limits_{-\pi }^{\pi }{\frac{\left| \sum\limits_{k=0}^{\infty }{}{{({{({{\mathbf{F}}_{\mu }})}^{-1}}{{[D_{N+\mu n}^{\mu }{{a}^{(2)}}]}_{\infty }})}_{k}}{{e}^{i\lambda k}} {\left. -\sum\limits_{k=-\infty }^{N+\mu n}{}\,{{({{({{\mathbf{G}}_{\mu }})}^{-1}}{{[D_{N+\mu n}^{\mu }{{a}^{(2)}}]}_{\infty }})}_{N+\mu n-k}}{{e}^{i\lambda k}} \right|}^{2} \right.}{{{\lambda }^{-2n}}|1-{{e}^{i\lambda \mu }}{{|}^{2n}}g(\lambda )}}\,d\lambda .
\end{multline}

Therefore, the following theorem holds.

\begin{thm} \label{theorem04}
 Let the uncorrelated stochastic sequences $\{\xi (m),m\in Z\}$ and $\{\eta (m),m\in Z\}$ determine the stationary $n$-th increments ${{\xi }^{(n)}}(m,\mu )$ and ${{\eta }^{(n)}}(m,\mu )$ with absolutely continuous spectral functions $F(\lambda )$ and $G(\lambda )$, which have spectral densities $f(\lambda )$ and $g(\lambda )$ satisfying conditions (15) and (16). Then the optimal linear estimate $\hat{A}_{N+1}^{N+\mu n}\xi $ of the value of the functional $A_{N+1}^{N+\mu n}\xi $ from the unknown values of $\xi (m)$, $m\in \{N+1,N+2,\ldots ,N+\mu n\}$, based on the known observations of the sequence $\xi (m)$ at points $m=-1,-2,\ldots $ and the sequence $\xi (m)+\eta (m)$ at points $m=N+1,N+2,\ldots $ is determined by formula (30). The spectral characteristic $\tilde{h}_{\mu }^{(a)}(\lambda )=(\tilde{h}_{\mu }^{(a,1)}(\lambda ),\tilde{h}_{\mu }^{(a,2)}(\lambda ))$ of the optimal estimate $\hat{A}_{N+1}^{N+\mu n}\xi $ is calculated by formulas (31) and (32). The magnitude of the mean square error $\Delta (f,g,\hat{A}_{N+1}^{N+\mu n})$ of the estimate $\hat{A}_{N+1}^{N+\mu n}\xi $ is calculated by formula (33).
\end{thm}

\section{{Minimax (robust) method of interpolation from noisy observations}}
The value of the mean square error $\Delta (h_{\mu }^{(a)}(f,g);f,g):=\Delta (f,g;{{\hat{A}}_{N}})$ and the spectral characteristic $h_{\mu }^{(a)}(\lambda )=(h_{\mu }^{(a,1)}(\lambda ),h_{\mu }^{(a,2)}(\lambda ))$ of the optimal linear estimate ${{\hat{A}}_{N}}\xi $ of the functional ${{A}_{N}}\xi $ from unknown values $\xi (m)$, $m=0,1,\ldots ,N$, based on observations of the sequence $\xi (k)$ at points $k\le -1$ and observations of the sequence $\xi (k)+\eta (k)$ at points $k\ge N+1$ are determined by formulas (25), (26) and (27),
if the spectral densities $f(\lambda )$ and $g(\lambda )$ of stochastic sequences $\xi (m)$ and $\eta (m)$ are known.
In the case where only a set $D={{D}_{f}}\times {{D}_{g}}$ of admissible densities is given, the minimax approach to the problem of estimating the functional is applied, i.e., an estimate is determined that minimizes the value of the mean square error simultaneously for all pairs of spectral densities from the class $D={{D}_{f}}\times {{D}_{g}}$.

\begin{ozn} \label{def3}
 For a given class of spectral densities $D={{D}_{f}}\times {{D}_{g}}$, the spectral densities ${{f}_{0}}(\lambda )\in {{D}_{f}}$, ${{g}_{0}}(\lambda )\in {{D}_{g}}$ are called the least favorable in the class $D$ for the optimal linear interpolation of the functional ${{A}_{N}}\xi $, if
	\[\Delta ({{f}_{0}},{{g}_{0}})=\Delta ({{h}_{\mu }}({{f}_{0}},{{g}_{0}});{{f}_{0}},{{g}_{0}})=\underset{(f,g)\in {{D}_{f}}\times {{D}_{g}}}{ {\max}}\,\Delta ({{h}_{\mu }}(f,g);f,g).\]
\end{ozn}

\begin{ozn} \label{def4} For a given class of spectral densities $D={{D}_{f}}\times {{D}_{g}}$, the spectral characteristic ${{h}^{0}}({{e}^{i\lambda }})$ of the optimal estimate of the functional $A_N\xi $ is called minimax (robust) if
	\[{{h}^{0}}({{e}^{i\lambda }})\in {{H}_{D}}=H_{D}^{(1)}\times H_{D}^{(2)}=\bigcap\limits_{(f,g)\in {{D}_{f}}\times {{D}_{g}}}{}\,L_{2}^{0-}(f)\times L_{2}^{(N+\mu n)+}(f+g),\]
	\[\underset{h\in {{H}_{D}}}{ {\min}}\,\underset{(f,g)\in {{D}_{f}}\times {{D}_{g}}}{ {\max}}\,\Delta ({{h}_{\mu }};f,g)=\underset{(f,g)\in {{D}_{f}}\times {{D}_{g}}}{ {\max}}\,\Delta ({{h}^{0}};f,g).\]
\end{ozn}

\begin{lem} \label{lem1}
The spectral densities $f^0(\lambda)\in{D}_f$ and $g^0(\lambda)\in{D}_g$ which satisfy  conditions \eqref{15}, \eqref{16},
are the  least favourable spectral densities in the class $D={{D}_{f}}\times {{D}_{g}}$ for the optimal
linear interpolation of the functional $A_N\xi$ based on observations
of the stochastic sequence $\xi(m)$ at points  $m=-1,-2,\ldots$ and observations of the stochastic sequence  $\xi(m)+\eta(m)$ at points  $m=N+1,N+2,\ldots$
if the matrices $ (\me F_{\mu})^0$, $(\me G_{\mu})^0$ defined by
the Fourier coefficients of the functions
\[
    \dfrac{\lambda^{2n}}{|1-e^{i\lambda\mu}|^{2n}f^0(\lambda)}, \quad
    \dfrac{\lambda^{2n}}{|1-e^{i\lambda\mu}|^{2n}g^0(\lambda)},\]
determine a solution of the constrained optimization problem
\begin{multline}\label{34}
    \max_{(f,g)\in \mathcal{D}_f\times\md D_g}
    \ld(\ip
    \frac{\lambda^{2n}\ld|
    C_{\mu}(e^{i\lambda })\rd|^2}{|1-e^{i\lambda\mu}|^{2n}f(\lambda)}d\lambda+\ip
    \frac{\lambda^{2n}\ld|
    C_{\mu}(e^{i\lambda })-E_{\mu}(e^{i\lambda})\rd|^2}{|1-e^{i\lambda\mu}|^{2n}g(\lambda)}d\lambda\rd)=
\\
={\int_{-\pi}^{\pi}}
    \frac{\lambda^{2n}\ld|
    C_{\mu}^0(e^{i\lambda })\rd|^2}{|1-e^{i\lambda\mu}|^{2n}f^0(\lambda)}d\lambda
    +\ip
    \frac{\ld|
    C_{\mu}^0(e^{i\lambda })-E_{\mu}^0(e^{i\lambda})\rd|^2}
    {\lambda^{-2n}|1-e^{i\lambda\mu}|^{2n}g^0(\lambda)}d\lambda,
    \end{multline}
\[
    C_{\mu}^0(e^{i\lambda })=\sum\limits_{k=0}^{\infty}
    \ld((\me F_{\mu}^{0})^{-1}\me [D_N^{\mu}\me     a_N]_{\infty}\rd)_ke^{i\lambda
    k},\]
\[
     E_{\mu}^0(e^{i\lambda})=\sum\limits_{k=-\infty}^{N+\mu n}\ld((\me G_{\mu}^{0})^{-1}\me [D_N^{\mu}\me a_N]_{\infty}\rd)_{N+\mu n-k}e^{i\lambda
k}.\]
The minimax spectral characteristic $h^0=h_{\mu}(f^0,g^0)$ is calculated by formulas \eqref{17}, \eqref{18} if
$h_{\mu}(f^0,g^0)\in H_{\mathcal{D}}$.
\end{lem}

For more detailed analysis of properties of the least favorable spectral densities and the minimax-robust spectral characteristics we observe that the minimax spectral characteristic $h^0$ and the  pair of the least favourable spectral densities $(f^0,g^0)$
form a saddle point of the functional $\Delta(h;f,g)$ on the set
$H_{\mathcal{D}}\times\mathcal{D}$.
The saddle point inequalities
\[\Delta (h;{{f}_{0}},{{g}_{0}})\ge \Delta ({{h}^{0}};{{f}_{0}},{{g}_{0}})\ge \Delta ({{h}^{0}};f,g)\quad \forall f\in {{D}_{f}},\forall g\in {{D}_{g}},\forall h\in {{H}_{D}}\]
are satisfied if ${{h}^{0}}={{h}_{\mu }}({{f}_{0}},{{g}_{0}})$ and ${{h}_{\mu }}({{f}_{0}},{{g}_{0}})\in {{H}_{D}}$, where $({{f}_{0}},{{g}_{0}})$ 
is a solution of the constrained optimization problem
\[\tilde{\Delta }(f,g)=-\Delta ({{h}_{\mu }}({{f}_{0}},{{g}_{0}});f,g)\to \inf,\quad (f,g)\in D,\]
	\[\Delta ({{h}_{\mu }}({{f}_{0}},{{g}_{0}});f,g)=\int\limits_{-\pi }^{\pi }{\frac{ {\left| \sum\limits_{k=0}^{\infty }{}\,{{({{(\mathbf{F}_{\mu }^{0})}^{-1}}{{[D_{N}^{\mu }{{a}^{(1)}}]}_{\infty }})}_{k}}{{e}^{i\lambda k}} \right|}^{2}}{{{\lambda }^{-2n}}|1-{{e}^{i\lambda \mu }}{{|}^{2n}}f_{0}^{2}(\lambda )}f(\lambda )d\lambda }\,+\]
	\[+\int\limits_{-\pi }^{\pi }{\frac{ {\left| \sum\limits_{k=0}^{\infty }{}\,{{({{(\mathbf{F}_{\mu }^{0})}^{-1}}{{[D_{N}^{\mu }{{a}^{(1)}}]}_{\infty }})}_{k}}{{e}^{i\lambda k}}-\sum\limits_{k=-\infty }^{N+\mu n}{}\,{{({{(\mathbf{G}_{\mu }^{0})}^{-1}}{{[D_{N}^{\mu }{{a}^{(1)}}]}_{\infty }})}_{N+\mu n-k}}{{e}^{i\lambda k}} \right|}^{2}}{{{\lambda }^{-2n}}|1-{{e}^{i\lambda \mu }}{{|}^{2n}}g_{0}^{2}(\lambda )}g(\lambda )d\lambda }\,.\]
The last problem is equivalent to the unconditional extremum problem
	\[{{\Delta }_{D}}(f,g)=\tilde{\Delta }(f,g)+\delta (f,g|{{D}_{f}}\times {{D}_{g}})\to inf,\]
This optimization problem is equivalent to the unconstrained optimization problem:
\[
    \Delta_{{D}}(f,g)=\widetilde{\Delta}(f,g)+ \delta(f,g|{D}_f\times
    {D}_g)\to\inf,\]
where $\delta(f,g|{D}_f\times
{D}_g)$ is the indicator function  of the set
${D}_f\times{D}_g$.
 A solution $(f^0,g^0)$  of the unconstrained optimization problem is characterized by the condition
  $0\in\partial\Delta_{{D}}(f^0,g^0)$ which is the necessary and sufficient  condition under which the pair  $(f^0,g^0)$ belongs to the set of minimums of the convex  functional $\Delta_{{D}}(f,g)$ [26], [27].  Here $\partial\Delta_{{D}}(f^0,g^0)$ denotes a subdifferential of the functional $\Delta_{{D}}(f,g)$ at the point $(f,g)=(f^0,g^0)$.

\section{{The least favorable densities in the class $D_{0,\mu }^{-}\times D_{0,\mu }^{-}$ }}

Consider the problem of the minimax (robust) interpolation of the functional $A_N\xi$  which depends on the unknown values of the sequence $\xi(m)$ with stationary increments
based on observations
of the stochastic sequence $\xi(m)$ at points  $m=-1,-2,\ldots$ and observations of the stochastic sequence  $\xi(m)+\eta(m)$ at points  $m=N+1,N+2,\ldots$,
under the condition that  the spectral densities are not known exactly while there is specified the set ${D}={D}_{0,\mu}^-\times{D}_{0,\mu}^-$ of admissible spectral densities, where
\[{{D}_{f}}=D_{0,\mu }^{-}=\left\{ f(\lambda )|\frac{1}{2\pi }\int\limits_{-\pi }^{\pi }{\,\frac{{{\lambda }^{2n}}}{|1-{{e}^{i\lambda \mu }}{{|}^{2n}}f(\lambda )}d\lambda }\ge {{P}_{1}} \right\},\,\,
\]
\[{{D}_{g}}=D_{0,\mu }^{-}=\left\{ g(\lambda )|\frac{1}{2\pi }\int\limits_{-\pi }^{\pi }{\frac{{{\lambda }^{2n}}}{|1-{{e}^{i\lambda \mu }}{{|}^{2n}}g(\lambda )}d\lambda }\,\ge {{P}_{2}} \right\}.\]

From the condition $0\in \partial {{\Delta }_{D}}({{f}_{0}},{{g}_{0}})$ we obtain the following relations that determine the least favorable spectral densities
\begin{equation}\label{35}
{\left| \sum\limits_{k=0}^{\infty }{}\,{{({{(\mathbf{F}_{\mu }^{0})}^{-1}}{{[D_{N}^{\mu }{{a}^{(1)}}]}_{\infty }})}_{k}}{{e}^{i\lambda k}} \right|}^{2}=p_{1}^{2},
\end{equation}
\begin{equation}\label{36}
{\left| \sum\limits_{k=0}^{\infty }{}\,{{({{(\mathbf{F}_{\mu }^{0})}^{-1}}{{[D_{N}^{\mu }{{a}^{(1)}}]}_{\infty }})}_{k}}{{e}^{i\lambda k}}-\sum\limits_{k=-\infty }^{N+\mu n}{}\,{{({{(\mathbf{G}_{\mu }^{0})}^{-1}}{{[D_{N}^{\mu }{{a}^{(1)}}]}_{\infty }})}_{N+\mu n-k}}{{e}^{i\lambda k}} \right|}^{2}=p_{2}^{2}.
\end{equation}

Let us denote by $\mathbf{p}_{1}^{\mu }$ the element of the space ${{\ell }_{2}}$ with coordinates $p_{1}^{\mu }(0)={{p}_{1}}$ and $p_{1}^{\mu }(k)=0$, $k>0$, and by $\mathbf{p}_{2}^{\mu }$ the element of the space ${{\ell }_{2}}$ with coordinates $p_{2}^{\mu }(N+\mu n)={{p}_{1}}-{{p}_{2}}$ and $p_{2}^{\mu }(k)=0$, $k\ne N+\mu n$. Then equations (35) satisfy the solutions of equation
	\[\mathbf{F}_{\mu }^{0}\mathbf{p}_{1}^{\mu }={{[D_{N}^{\mu }{{a}^{(1)}}]}_{\infty }},\]
and equations (36) satisfy the solutions of the equation
	\[\mathbf{G}_{\mu }^{0}\mathbf{p}_{2}^{\mu }={{[D_{N}^{\mu }{{a}^{(1)}}]}_{\infty }}.\]
The last two equations are equivalent to the following
\begin{equation}\label{37}
\left( {{(\mathbf{G}_{\mu }^{c})}^{0}}+{{(\mathbf{F}_{\mu }^{c})}^{0}} \right){{({{(\mathbf{G}_{\mu }^{e})}^{0}})}^{-1}}{{(\mathbf{G}_{\mu }^{c})}^{0}}\mathbf{p}_{1}^{\mu }={{(\mathbf{G}_{\mu }^{e})}^{0}}\mathbf{p}_{1}^{\mu }+{{[D_{N}^{\mu }{{a}^{(1)}}]}_{\infty }},
\end{equation}
\begin{equation}\label{38}
\left( {{(\mathbf{G}_{\mu }^{c})}^{0}}+{{(\mathbf{F}_{\mu }^{c})}^{0}}-{{(\mathbf{G}_{\mu }^{e})}^{0}}{{({{(\mathbf{G}_{\mu }^{c})}^{0}})}^{-1}}{{(\mathbf{G}_{\mu }^{e})}^{0}} \right){{[D_{N}^{\mu }{{a}^{(1)}}]}_{\infty }}=\mathbf{p}_{2}^{\mu },
\end{equation}
where the matrices ${{(\mathbf{F}_{\mu }^{c})}^{0}}$, ${{(\mathbf{F}_{\mu }^{e})}^{0}}$, ${{(\mathbf{G}_{\mu }^{c})}^{0}}$, ${{(\mathbf{G}_{\mu }^{e})}^{0}}$
constructed using Fourier coefficients
	\[f_{\mu }^{0}(k)=\frac{1}{2\pi }\int\limits_{-\pi }^{\pi }{\frac{{{\lambda }^{2n}}{{e}^{-i\lambda k}}}{|1-{{e}^{i\lambda \mu }}{{|}^{2n}}{{f}^{0}}(\lambda )}d\lambda }\,,\quad g_{\mu }^{0}(k)=\frac{1}{2\pi }\int\limits_{-\pi }^{\pi }{\,\frac{{{\lambda }^{2n}}{{e}^{-i\lambda k}}}{|1-{{e}^{i\lambda \mu }}{{|}^{2n}}{{g}^{0}}(\lambda )}d\lambda },\quad k\in Z.\]
of the functions $$\frac{{{\lambda }^{2n}}}{|1-{{e}^{i\lambda \mu }}{{|}^{2n}}{{f}^{0}}(\lambda )},\quad \frac{{{\lambda }^{2n}}}{|1-{{e}^{i\lambda \mu }}{{|}^{2n}}{{g}^{0}}(\lambda )}.$$
The coefficients ${{p}_{1}}$ and ${{p}_{2}}$ are determined by the conditions
$$\frac{1}{2\pi }\int\limits_{-\pi }^{\pi }{\frac{{{\lambda }^{2n}}{{e}^{-i\lambda k}}}{|1-{{e}^{i\lambda \mu }}{{|}^{2n}}{{f}^{0}}(\lambda )}d\lambda }={{P}_{1}},\quad \frac{1}{2\pi }\int\limits_{-\pi }^{\pi }{\frac{{{\lambda }^{2n}}{{e}^{-i\lambda k}}}{|1-{{e}^{i\lambda \mu }}{{|}^{2n}}{{g}^{0}}(\lambda )}d\lambda }={{P}_{2}}.$$

Let the sequences $\{f_{\mu }^{0}(k):k\ge 0\}$, $\{g_{\mu }^{0}(k):k\ge 0\}$ satisfy equations (37), (38). Then the least favorable densities ${{f}^{0}}(\lambda )$ and ${{g}^{0}}(\lambda )$ are represented in the form
\begin{equation}\label{39}
{{f}^{0}}(\lambda )={{\lambda }^{2n}}{\left( |1-{{e}^{i\lambda \mu }}{{|}^{2n}}\sum\limits_{k=-\infty }^{\infty }{}\,f_{\mu }^{0}(|k|){{e}^{i\lambda k}} \right)}^{-1}={\left| \frac{{{(i\lambda )}^{n}}}{{{(1-{{e}^{-i\lambda \mu }})}^{n}}\sum\limits_{k=0}^{\infty }{}\,{{\gamma }_{\mu }}(k){{e}^{-i\lambda k}}} \right|}^{2},
\end{equation}
\begin{equation}\label{40}
{{g}^{0}}(\lambda )={{\lambda }^{2n}}{\left( |1-{{e}^{i\lambda \mu }}{{|}^{2n}}\sum\limits_{k=-\infty }^{\infty }{}\,g_{\mu }^{0}(|k|){{e}^{i\lambda k}} \right)}^{-1}={\left| \frac{{{(i\lambda )}^{n}}}{{{(1-{{e}^{-i\lambda \mu }})}^{n}}\sum\limits_{k=0}^{\infty }{}\,{{\zeta }_{\mu }}(k){{e}^{-i\lambda k}}} \right|}^{2}.
\end{equation}

\begin{thm} \label{theorem05}
Let the coefficients $\{f_{\mu }^{0}(k):k\ge 0\}$, $\{g_{\mu }^{0}(k):k\ge 0\}$ satisfy the system (37), (38). The least favorable densities in the class $D={{D}_{f}}\times {{D}_{g}}$ for constructing a linear estimate of the functional ${{A}_{N}}\xi $ from observations of the sequence $\xi (k)$ at points $k\le -1$ and bservations of the sequence $\xi (k)+\eta (k)$ at points $k\ge N+1$, where $\eta (k)$ is an uncorrelated sequence with $\xi (k)$, have the form (39), (40). The minimax estimate of the functional ${{A}_{N}}\xi $ is calculated by formula (24) with $h_{\mu }^{(a)}(\lambda )={{h}^{0}}(\lambda )$. The minimax spectral characteristic $h_{\mu }^{0}(\lambda )={{h}_{\mu }}({{f}_{0}},g)$ is calculated by formulas (25), (26).
\end{thm}

\begin{thm} \label{theorem06}
Let the spectral density $g(\lambda )$ be known, and the coefficients $\{f_{\mu }^{0}(k):k\ge 0\}$ satisfy the equation
\begin{equation}\label{41}
{{(\mathbf{F}_{\mu }^{c})}^{0}}{{(\mathbf{G}_{\mu }^{e})}^{-1}}\mathbf{G}_{\mu }^{c}\mathbf{p}_{1}^{\mu }=\mathbf{G}_{\mu }^{e}\mathbf{p}_{1}^{\mu }-\mathbf{G}_{\mu }^{c}{{(\mathbf{G}_{\mu }^{e})}^{-1}}\mathbf{G}_{\mu }^{c}\mathbf{p}_{1}^{\mu }+{{[D_{N}^{u}{{a}^{(1)}}]}_{\infty }},
\end{equation}
The least favorable density in the class $D={{D}_{f}}$ for constructing a linear estimate of the functional ${{A}_{N}}\xi $ from observations of the sequence $\xi (k)$ at points $k\le -1$ and
 bservations of the sequence $\xi (k)+\eta (k)$ at points $k\ge N+1$, where $\eta (k)$ is an uncorrelated sequence with $\xi (k)$, has the form (39).
 The minimax estimate of the functional ${{A}_{N}}\xi $ is calculated by formula (24) with $h_{\mu }^{(a)}(\lambda )={{h}^{0}}(\lambda )$.
 The minimax spectral characteristic ${{h}^{0}}(\lambda )={{h}_{\mu }}({{f}_{0}},g)$ is calculated by formulas (25), (26).
\end{thm}

Suppose that the $n$-th increment ${{\eta }^{(n)}}(k,\mu )$ of the sequence $\eta (k)$ is a sequence of independent identically distributed random variables with mean 0 and variance ${{\sigma }^{2}}$ (white noise). In this case, the matrix $\mathbf{G}_{\mu }^{e}={{\sigma }^{2}}I$, where $I$ is the unit operator, and the matrix $\mathbf{G}_{\mu }^{c}$ is given by the elements ${{(\mathbf{G}_{\mu }^{c})}_{k,l}}={{\sigma }^{2}}$ for $k+l=N+\mu n$, $k,l\ge 0$, and ${{(\mathbf{G}_{\mu }^{c})}_{k,l}}=0$ otherwise. Equation (47) is equivalent to the following system
	\[{{p}_{1}}f_{\mu }^{0}(k)={{({{[D_{N}^{\mu }{{a}^{(1)}}]}_{\infty }})}_{k}},\quad k\ge 0.\]
From the condition
$$\frac{1}{2\pi }\int\limits_{-\pi }^{\pi }{\frac{{{\lambda }^{2n}}{{e}^{-i\lambda k}}}{|1-{{e}^{i\lambda \mu }}{{|}^{2n}}{{f}^{0}}(\lambda )}d\lambda }={{P}_{1}}$$
we get  ${{p}_{1}}=P_{1}^{-1}{{({{[D_{N}^{\mu }{{a}^{(1)}}]}_{\infty }})}_{0}}$, $f_{\mu }^{0}(k)={{P}_{1}}{{(D_{N}^{\mu }{{a}^{(1)}})}_{k}}/{{(D_{N}^{\mu }{{a}^{(1)}})}_{0}}$ for $k=0,1,\ldots ,N$ and $f_{\mu }^{0}(k)=0$ for $k>N$.

Therefore, if the sequence $({{(D_{N}^{\mu }{{a}^{(1)}})}_{0}},{{(D_{N}^{\mu }{{a}^{(1)}})}_{1}},\ldots ,{{(D_{N}^{\mu }{{a}^{(1)}})}_{N}})$ is positive, the least favorable density can be calculated by the formula
\begin{equation}\label{42}
{{f}^{0}}(\lambda )={{\lambda }^{2n}}{\left( |1-{{e}^{i\lambda \mu }}{{|}^{2n}}\sum\limits_{k=-N}^{N}{}\,f_{\mu }^{0}(|k|){{e}^{i\lambda k}} \right)}^{-1}={\left| \frac{{{(i\lambda )}^{n}}}{{{(1-{{e}^{-i\lambda \mu }})}^{n}}\sum\limits_{k=0}^{N}{}\,{{\gamma }_{\mu }}(k){{e}^{-i\lambda k}}} \right|}^{2}.
\end{equation}

\begin{nas} \label{corr2}
Let the $n$-th increment of the sequence $\eta (k)$ be white noise, and the sequence $({{(D_{N}^{\mu }{{a}^{(1)}})}_{0}},{{(D_{N}^{\mu }{{a}^{(1)}})}_{1}},\ldots ,{{(D_{N}^{\mu }{{a}^{(1)}})}_{N}})$ be positive. Then the least favorable density in the class $D={{D}_{f}}$ for constructing a linear estimate of the functional ${{A}_{N}}\xi $ from observations of the sequence $\xi (k)$ at points $k\le -1$ and the sequence $\xi (k)+\eta (k)$ at points $k\ge N+1$ has the form (42), where the coefficients $f_{\mu }^{0}(k)={{P}_{1}}{{(D_{N}^{\mu }{{a}^{(1)}})}_{k}}/{{(D_{N}^{\mu }{{a}^{(1)}})}_{0}}$ for $k=0,1,\ldots ,N$.
\end{nas}

\section{{The least favorable densities in the class $D_{M,\mu }^{-}\times D_{M,\mu }^{-}.$ }}

Consider the problem of the minimax (robust) interpolation of the functional $A_N\xi$  which depends on the unknown values of the sequence $\xi(m)$ with stationary increments
based on observations
of the stochastic sequence $\xi(m)$ at points  $m=-1,-2,\ldots$ and observations of the stochastic sequence  $\xi(m)+\eta(m)$ at points  $m=N+1,N+2,\ldots$,
under the condition that  the spectral densities are
is not known exactly while there is specified the set $D={{D}_{f}}\times {{D}_{g}}=D_{M,\mu }^{-}\times D_{M,\mu }^{-}$ of admissible  spectral densities, where
	\[{{D}_{f}}=D_{M,n}^{-}=\left\{ f(\lambda )|\frac{1}{2\pi }\int\limits_{-\pi }^{\pi }{\frac{{{\lambda }^{2n}}}{|1-{{e}^{i\lambda \mu }}{{|}^{2n}}f(\lambda )}\cos(m\lambda )d\lambda }\,={{r}_{1}}(m),m=0,1.\ldots ,M \right\},\]
	\[{{D}_{g}}=D_{M,n}^{-}=\left\{ g(\lambda )|\frac{1}{2\pi }\int\limits_{-\pi }^{\pi }{\frac{{{\lambda }^{2n}}}{|1-{{e}^{i\lambda \mu }}{{|}^{2n}}g(\lambda )}\cos(m\lambda )d\lambda }\,={{r}_{2}}(m),m=0,1.\ldots ,M \right\},\]
where $\{{{r}_{1}}(m),m=0,1,\ldots ,M\}$ and $\{{{r}_{2}}(m),m=0,1,\ldots ,M\}$ are strictly positive sequences. From the condition $0\in \partial {{\Delta }_{D}}({{f}^{0}})$ we find the equation
\begin{equation}\label{43}
{\left| \sum\limits_{k=0}^{\infty }{}\,{{({{(\mathbf{F}_{\mu }^{0})}^{-1}}{{[D_{N}^{\mu }{{a}^{(1)}}]}_{\infty }})}_{k}}{{e}^{i\lambda k}} \right|}^{2}=\sum\limits_{m=0}^{M}{}\,{{\alpha }_{1,m}}\cos(m\lambda) ={\left| \sum\limits_{m=0}^{M}{}\,{{p}_{1}}(m){{e}^{i\lambda m}} \right|}^{2},
\end{equation}
\begin{multline}\label{44}
{\left| \sum\limits_{k=0}^{\infty }{}\,{{({{(\mathbf{F}_{\mu }^{0})}^{-1}}{{[D_{N}^{\mu }{{a}^{(1)}}]}_{\infty }})}_{k}}{{e}^{i\lambda k}}-\sum\limits_{k=-\infty }^{N+\mu n}{}\,{{({{(\mathbf{G}_{\mu }^{0})}^{-1}}{{[D_{N}^{\mu }{{a}^{(1)}}]}_{\infty }})}_{N+\mu n-k}}{{e}^{i\lambda k}} \right|}^{2}=\\
=\sum\limits_{m=0}^{M}{}\,{{\alpha }_{2,m}}\cos(m\lambda) =
{\left| \sum\limits_{m=0}^{M}{}\,{{p}_{2}}(m){{e}^{i\lambda m}} \right|}^{2}.
\end{multline}
where ${{\alpha }_{j,m}}$, $j=1.2$, $m=0.1,\ldots ,M$, are Lagrange multipliers. We denote by $\mathbf{p}_{M}^{1,\mu }$ an element of the space ${{\ell }_{2}}$ with coordinates $p_{M}^{1,\mu }(k)={{p}_{1}}(k)$, $k=0,1,\ldots ,M$, and $p_{M}^{1,\mu }(k)=0$, $k\ge M+1$. We denote by $\mathbf{p}_{M}^{2,\mu }$ another element of the space ${{\ell }_{2}}$. Consider two cases: $M\le N+\mu n$ and $M>N+\mu n$. For $M\le N+\mu n$, we set $p_{M}^{2,\mu }(k)={{p}_{1}}(M-k)-{{p}_{2}}(M-k)$, $k=0,1,\ldots,M$, and $p_{M}^{2,\mu }(k)=0$, $k>M$. If $M>N+\mu n$, we set $p_{M}^{2,\mu }(k)={{p}_{1}}(N+\mu n-k)-{{p}_{2}}(N+\mu n-k)$, $k=0,1,\ldots ,N+\mu n$, and $p_{M}^{2,\mu }(k)=0$, $k>N+\mu n$. Let also in the case $M>N+\mu n$ the equalities ${{p}_{1}}(k)={{p}_{2}}(k)$, $k=N+\mu n,\ldots ,M$ hold. Then equations (43) and (44) satisfy the solutions of equations
	\[\mathbf{F}_{\mu }^{0}\mathbf{p}_{M}^{1,\mu }={{[D_{N}^{\mu }{{a}^{(1)}}]}_{\infty }},\quad \mathbf{G}_{\mu }^{0}\mathbf{p}_{M}^{2,\mu }={{[D_{N}^{\mu }{{a}^{(1)}}]}_{\infty }}.\]
The last two equations are equivalent to the following
\begin{equation}\label{45}
\left( {{(\mathbf{G}_{\mu }^{c})}^{0}}+{{(\mathbf{F}_{\mu }^{c})}^{0}} \right){{({{(\mathbf{G}_{\mu }^{e})}^{0}})}^{-1}}{{(\mathbf{G}_{\mu }^{c})}^{0}}\mathbf{p}_{M}^{1,\mu }={{(\mathbf{G}_{\mu }^{e})}^{0}}\mathbf{p}_{M}^{1,\mu }+{{[D_{N}^{\mu }{{a}^{(1)}}]}_{\infty }},
\end{equation}
\begin{equation}\label{46}
\left( {{(\mathbf{G}_{\mu }^{c})}^{0}}+{{(\mathbf{F}_{\mu }^{c})}^{0}}-{{(\mathbf{G}_{\mu }^{e})}^{0}}{{({{(\mathbf{G}_{\mu }^{c})}^{0}})}^{-1}}{{(\mathbf{G}_{\mu }^{e})}^{0}} \right){{[D_{N}^{\mu }{{a}^{(1)}}]}_{\infty }}=\mathbf{p}_{M}^{2,\mu }.
\end{equation}
To determine the coefficients ${{p}_{1}}(k)$ and ${{p}_{2}}(k)$, $k=0,1,\ldots ,M$, we use the equalities $f_{\mu }^{0}(m)={{r}_{1}}(m)$, $g_{\mu }^{0}(m)={{r}_{2}}(m)$, $m=0,1,\ldots ,M$, and in the case $M>N+\mu n$ we must additionally require ${{p}_{1}}(k)={{p}_{2}}(k)$ for $k=N+\mu n,\ldots ,M$.

Let the sequences $\{f_{\mu }^{0}(k):k\ge 0\}$, $\{g_{\mu }^{0}(k):k\ge 0\}$ satisfy equations (45), (46).

Then the least favorable densities ${{f}^{0}}(\lambda )$ and ${{g}^{0}}(\lambda )$ can be represented as (39) and (40), respectively.

\begin{thm} \label{theorem07}
  Let the coefficients $\{f_{\mu }^{0}(k):k\ge 0\}$, $\{g_{\mu }^{0}(k):k\ge 0\}$ satisfy the system (45), (46). The least favorable densities in the class $D=D_{M,\mu }^{-}\times D_{M,\mu }^{-}$ for constructing a linear estimate of the functional ${{A}_{N}}\xi $ from observations of the sequence $\xi (k)$ at points $k\le -1$ and the sequence $\xi (k)+\eta (k)$ at points $k\ge N+1$, where $\eta (k)$ is a sequence uncorrelated with $\xi (k)$, have the form (39), (40).
  The minimax estimate of the functional ${{A}_{N}}\xi $ is calculated by formula (24) where $h_{\mu }^{(a)}(\lambda )={{h}^{0}}(\lambda )$.
  The minimax spectral characteristic $h_{\mu }^{0}(\lambda )={{h}_{\mu }}({{f}_{0}},g)$ is calculated by formulas (25), (26).
\end{thm}
\begin{thm} \label{theorem08}
 Let the spectral density $g(\lambda )$ be known, and the coefficients $\{f_{\mu }^{0}(k):k\ge 0\}$ satisfy the equation
\begin{equation}\label{47}
{{(\mathbf{F}_{\mu }^{c})}^{0}}{{(\mathbf{G}_{\mu }^{e})}^{-1}}\mathbf{G}_{\mu }^{c}\mathbf{p}_{M}^{1,\mu }=\mathbf{G}_{\mu }^{e}\mathbf{p}_{M}^{1,\mu }-\mathbf{G}_{\mu }^{c}{{(\mathbf{G}_{\mu }^{e})}^{-1}}\mathbf{G}_{\mu }^{c}\mathbf{p}_{M}^{1,u}+{{[D_{N}^{\mu }{{a}^{(1)}}]}_{\infty }},
\end{equation}
The least favorable density in the class $D=D_{M,\mu }^{-}$ for constructing a linear estimate of the functional ${{A}_{N}}\xi $ from observations of the sequence $\xi (k)$ at points $k\le -1$ and the sequence $\xi (k)+\eta (k)$ at points $k\ge N+1$, where $\eta (k)$ is an uncorrelated sequence with $\xi (k)$, has the form (39).
The minimax estimate of the functional ${{A}_{N}}\xi $ is calculated by formula (24) with $h_{\mu }^{(a)}(\lambda )={{h}^{0}}(\lambda )$.
The minimax spectral characteristic ${{h}^{0}}(\lambda )={{h}_{\mu }}({{f}_{0}},g)$ is calculated by formulas (25), (26).
\end{thm}

\section{{Conclusions}}
 The paper provides formulas for calculating the spectral characteristics and the mean square errors of optimal estimates of the functionals $A_{N}^{{}}\xi =\sum\limits_{k=0}^{N}{a(k)\xi (k)}$ and $A_{N+1}^{N+\mu n}\xi =\sum\limits_{k=N+1}^{N+\mu n}{a(k)\xi (k)}$ from unknown values of the stochastic sequence $\xi (m)$ with stationary $n$-th increments based on observations of the sequence $\xi (k)$ at points $k=-1,-2,\ldots $ and observations of the sequence $\xi (k)+\eta (k)$ at points $k=N+1,N+2,\ldots $, where $\eta (k)$ is a stochastic sequence with stationary $n$ increments, uncorrelated with the sequence $\xi (k)$, in the case where the spectral densities of the sequences $\xi (k)$ and $\eta (k)$ are known.
 In the case where the spectral densities are unknown, while some sets of admissible spectral densities are given, the minimax estimation method was applied to solve the problem of interpolation of the functional $A_{N}^{{}}\xi$. For certain classes of admissible sets of spectral densities, relations were found that determine the least favorable spectral densities and the minimax spectral characteristics of the optimal estimate of the functional ${{A}_{N}}\xi $.


\begin{thebibliography}{99}

\bibitem{Yaglom:1955}
Yaglom, A. M. (1955). Correlation theory of stationary and related random processes with stationary nth
increments. Mat. Sbornik 37(1), 141-196.

\bibitem{Pinsk:1955}
Pinsker, M. S. (1955). The theory of curves with nth stationary increments in Hilbert spaces. Izvestiya
Akademii Nauk SSSR. Ser. Mat. 19(5), 319-344.

\bibitem{Pinsk:1954}
Pinsker, M. S.,  Yaglom, A. M. (1954). On linear extrapolation of random processes with nth stationary
increments. Doklady Akademii Nauk SSSR 94, 385-388.

\bibitem{Yaglom:1957}
Yaglom, A. M. (1957). Some classes of random fields in n-dimensional space related with random stationary
processes. Teor. Veroyatn. Primen. 2, 292-338.

\bibitem{Kolmogorov}
Kolmogorov, A.N. (1992). Selected works of A. N. Kolmogorov. Volume II: Probability theory and mathematical statistics.
Edited by A. N. Shiryayev. Transl. from the Russian by G. Lindquist.
Dordrecht etc.: Kluwer Academic Publishers, 584.

\bibitem{Wiener}
Wiener, N. (1966). Extrapolation, interpolation, and smoothing of
stationary time series. With engineering applications. Cambridge, Mass.: The M. I. T. Press, Massachusetts Institute of Technology, 163.

\bibitem{Yaglom:1987a}
Yaglom, A. M. (1987). Correlation theory of stationary and related random functions. Vol. 1: Basic results.
Springer Series in Statistics, Springer-Verlag, New York etc., 526.

\bibitem{Yaglom:1987b}
Yaglom, A. M. (1987). Correlation theory of stationary and related random functions. Vol. 2: Supplementary
notes and references. Springer Series in Statistics, Springer-Verlag, New York etc., 258.

\bibitem{Luz:2011}
Luz, M. (2011). Wiener-Kolmogorov predicting for stationary processes,
Visn., Mat. Mekh., Kyiv. Univ. Im. Tarasa Shevchenka no. 25, 26--29

\bibitem{Grenander}
Grenander, U. (1957). A prediction problem in game theory. Arkiv f\"or Matematik, 3, 371-379.

\bibitem{Franke:1985}
Franke, J. (1985). Minimax robust prediction of discrete time series. Z.
Wahrsch. Verw. Gebiete, 68, 337-364.

\bibitem{FrankePoor}
Franke, J.,   Poor, H. V. (1984). Minimax-robust filtering
and finite-length robust predictors. Robust and Nonlinear Time
Series Analysis. Lecture Notes in Statistics, Springer-Verlag
26, 87-126.

\bibitem{KassamPoor}
Kassam, S. A., Poor, H. V. (1985). Robust techniques for signal processing: A
survey. Proceedings of the IEEE, 73(3), 433-481.

\bibitem{Moklyachuk1991}
Moklyachuk M.P., (1991). Minimax filtration of linear transformations of stationary sequences. Ukr. Math. J.  43(1),  75-81.

\bibitem{Moklyachuk1984}
Moklyachuk, M. P., (1994). Stochastic autoregressive sequences and minimax interpolation. Theory Probability and Mathematical Statistics 48, 95-103.

\bibitem{Moklyachuk2000}
Moklyachuk, M. P., (2000). Robust procedures in time series analysis. Theory Stochastic Processes 6(3-4), 127-147.

\bibitem{Moklyachuk2001}
Moklyachuk, M. P., (2001). Game theory and convex optimization methods in robust estimation problems. Theory Stochastic Processes 7(1-2), 253-264.

\bibitem{Moklyachuk2008}
Moklyachuk, M. P., (2008). Robust estimations of functionals of stochastic processes. Vydavnycho-Poligrafichny\u\i\ Tsentr, Ky{\"\i}vsky\u\i\ Universytet, Ky{\"\i}v, 320.

\bibitem{Rozanov}
Rozanov, Yu. A. (1967). Stationary stochastic processes. Holden-Day, San Francisco, 211.

\bibitem{Moklyachuk:2012}
Moklyachuk, M. P.,   Masyutka, O. Yu. (2012). Minimax-robust
estimation technique for stationary stochastic processes, LAP
Lambert Academic Publishing, 296.

\bibitem{Dubovetska4}
Dubovets'ka, I. I.,  Moklyachuk, M. P. (2013).
Filtration of
linear functionals of periodically correlated sequences, Theor.
Theory of Probability and Mathematical Statistics, 86, 51-64.

\bibitem{Dubovetska1}
Dubovets'ka, I. I., Masyutka, O. Yu.,  Moklyachuk, M. P. (2012).
Interpolation of periodically correlated stochastic sequences.
Theory of Probability and Mathematical Statistics, 84, 43-56.

\bibitem{Luz2}
Luz, M. M.,    Moklyachuk, M. P. (2013). Interpolation of functionals of stochastic sequences with stationary
increments. Theor. Probability and  Math. Stat. 87, 117-133.

\bibitem{Karhunen}
Karhunen, K. (1947). Uber lineare Methoden in der Wahrscheinlichkeitsrechnung. Annales Academiae
Scientiarum Fennicae. Series A I. Mathematica 37, 3-79.

\bibitem{Gikhman:Skorokhod}
 Gikhman, I. I.,  Skorokhod, A. V. (2004). The theory of stochastic processes. I., Springer, Berlin.

\bibitem{Moklyachuk:2008nonsm}
Moklyachuk, M. P., (2008). Nonsmooth analysis and optimization.
Vydavnycho-Poligrafichny\u\i\ Tsentr, Ky{\"\i}vsky\u\i\ Universytet, Ky{\"\i}v, 400.

\bibitem{Pshenychn}
 Pshenichnyi, B. N. (1971). Necessary conditions for an extremum.
Pure and Applied mathematics. 4.; New York: Marcel Dekker, Inc. XVIII, 230.


\end{thebibliography}
\end{document}